\documentclass[a4paper,12pt]{article}

\usepackage{amsmath}
\usepackage{amsfonts} 
\usepackage{amssymb,amsthm}
\usepackage{color}
\definecolor{rouge}{rgb}{1,0,0}
\definecolor{gris}{gray}{0.5}

\usepackage{hyperref}

\usepackage{pb-diagram}  
\usepackage[cmtip,arrow]{xy}
\usepackage{pb-xy} 
\usepackage{enumerate}   

\usepackage{bbold} 

\usepackage[english]{babel}   
\usepackage{fullpage}       

\usepackage[utf8]{inputenc} 
\usepackage{lmodern}

\newtheorem{de}{Definition}[section]
\newtheorem{theo}[de]{Theorem}
\newtheorem{prop}[de]{Proposition}
\newtheorem{cor}[de]{Corollary}
\newtheorem{lem}[de]{Lemma}
\newtheorem{rem}[de]{Reminder}
\newtheorem{propdef}[de]{Proposition/Definition}

\numberwithin{equation}{section}

\newcommand{\bb}[1][R]{\mathbb{#1}}

\newcommand{\sbb}[1][1|1]{\mathbb{R}^{#1}}
\newcommand{\sbbpq}{\mathbb{R}^{p|q}}

\newcommand{\ca}[1][M]{{\mathcal{#1}}} 
\newcommand{\can}[1][N]{{\mathcal{#1}}} 
\newcommand{\ab}{{\alpha\beta}}
\newcommand{\al}{_\alpha}
\newcommand{\be}{_\beta}
\newcommand{\om}[1][M]{{\ca[O]_{\ca[#1]}}}
\newcommand{\on}[1][N]{{\ca[O]_{\ca[#1]}}}
\newcommand{\tm}[1][M]{\ca[T]_{\ca[#1]}}
\newcommand{\tn}[1][N]{\ca[T]_{\ca[#1]}}
\newcommand{\pa}[1][q_i]{\partial_{#1}}
\newcommand{\sig}[1][|q_i|]{(-1)^{#1}}
\newcommand{\pai}[1][i]{\partial_{q_#1}}
\newcommand{\paj}[1][j]{\partial_{q_#1}}
\newcommand{\pak}[1][k]{\partial_{q_#1}}
\newcommand{\paq}[1][i]{\partial_{q_#1}}
\newcommand{\tphi}[1][\Phi]{\ca[T]^{#1}}
\newcommand{\tphiup}[1][\Phi]{\ca[T]^{#1}}
\newcommand{\cat}{\ca=(M,\om)}
\newcommand{\cant}{\can=(N,\on)}

\begin{document}

\title{The geodesic flow on a Riemannian supermanifold}
\author{
{\bf St\'ephane Garnier\footnote{partially supported by the Universit\'e 
Franco-Allemande (DFH-UFA) via the IRTG 1133 ``Geometry and Analysis of Symmetries"} 
}\\
{\bf Tilmann Wurzbacher\footnote{Corresponding author, e-mail: \tt{wurzbacher@math.univ-metz.fr}}
\footnote{On leave of absence; present address: Fakult\"at f\"ur Mathematik, Ruhr-Universit\"at Bochum,
Germany}
}\\
     {\small Laboratoire de Math\'ematiques et Applications de Metz}\\
             {\small Universit\'e Paul Verlaine-Metz et C.\,N.\,R.\,S.}\\
                 {\small Ile du Saulcy, F-57045 Metz, France}\\
}
\date{July 9, 2011}
\maketitle

\centerline{  {\bf MSC (2010)} Primary: 58A50  Secondary: 53C22  $\cdot$ 53D25  $\cdot$  58C50 }

\begin{abstract}
We give a natural definition of geodesics on a Riemannian supermanifold $(\ca,g)$ and  extend the usual
geodesic flow on $T^*M$ associated to the underlying Riemannian manifold 
$(M,\widetilde{g})$ to a geodesic ``superflow" on $T^*\ca$. Integral curves of this flow  turn out 
to be in natural bijection with geodesics
on $\ca$. We also construct the corresponding exponential map and generalize 
the well-known faithful linearization of isometries to Riemannian supermanifolds.
\end{abstract}

\section{Introduction}
The geometry of supermanifolds is, after a first, pioneering, ``Storm and Stress"-period in the 70s and 80s, still in its 
youth but with by now many mathematicians studying in detail the fundaments as well as the applications of this theory.
Our work on Riemannian supermanifolds is motivated by at least two developements originating in physics and relying on Riemannian geometry
extended to the category of supermanifolds. First, we are interested in those aspects of the ``Stolz-Teichner program"
on the geometrisation of elliptic cohomology that relate to cobordisms of Riemannian supermanifolds and to supermanifolds in general (see \cite{Stolz-Teichner} for a detailed account). Of course, extending certain results of this article to super loop spaces (compare  \cite{Wurzbacher} for the ungraded case) would give a direct connection to the physics literature underlying elliptic cohomology, where the role of loop spaces is central (see, e.g., \cite{Witten}). Second, work of Zirnbauer and collaborators (see, e.g.,
\cite{Zirnbauer}) shows the surprising importance of the notion of ``supersymmetric spaces" for physics, notably
mesoscopic physics. In this context it is very natural to ask for a general theory of Riemannian supermanifolds with these spaces replacing the symmetric spaces of ordinary Riemannian geometry. The subject of supersymmetric spaces was, in fact, the main motivation of the work of Goertsches
(\cite{Goertsches}) on Riemannian supergeometry.

\vspace*{5mm}
The first cornerstone of Riemannian geometry is -after the metric itself- the notion of a geodesic curve. Upon extending 
these notions, metrics and geodesics,
to supermanifolds it looks highly reasonable to aim at fulfilling  the following three conditions:
\begin{enumerate}[(i)]
 \item geodesic supercurves should reduce to the usual geodesics in the case of classical (ungraded) Riemannian
 manifolds
 \item an isometry fixing a point and having a trivial linearization in this point should be the identity morphism
 \item geodesics should be closely related to integral curves of a ``geodesic flow" on the
 (co-)tangent bundle of the given Riemannian manifold.
  \end{enumerate}

The definition of a Riemannian metric on a Riemannian metric is fairly straight forward and seems to be generally 
accepted in the literature. (See Definition 2.7 below.)    
The notion of geodesic supercurves proposed by Goertsches fulfils (i) and (ii), but unfortunately we could
not establish (iii) in this approach. Since the
geodesic flow is a tool of great importance in the geometry and analysis on a Riemannian manifold, we found
it crucial not to give up condition (iii). Furthermore, the discussion of initial conditions in \cite{Goertsches}
(see the beginning of Subsection 4.4 there) does not seem to be fully conclusive. Let us underline that the equations chosen 
for defining geodesics in this article translate in local coordinates to a system of second-order ODEs, whereas
the geodesics of Goertsches give rise to 2nd-order equations and 1st-order equations (for the even resp. odd
coordinates). This explains, at least on an intuitive level, why there might be no geodesic flow related to
Goertsches' notion of geodesics.
\vspace*{5mm}

Let us now proceed to a short description of the contents of this article. 

\vspace*{3mm}
First of all, we underline that we work throughout in the sheaf-theoretic framework, i.e., in the approach to 
supermanifolds going back to
Kostant and the russian school around Berezin, Leites and Manin.

\vspace*{3mm}
We define, in the second section, geodesics as ``supercurves" that are autoparallel with respect to 
the natural even time direction:

\[
0= \displaystyle\frac{\nabla}{dt}\big(T\Phi^*(\pa[t])\big)=(\Phi^*\nabla)_{\pa[t]}(\pa[t] \circ \Phi^*),
\]

where $\Phi^*\nabla$ is the pull-back of the Levi-Civita connection, $\pa[t]$ is $\frac{\partial}{\partial t}$
and   $\pa[t]\circ \Phi^*$ is a canonically defined vector field along $\Phi$ (see the main text for more details).
The preceding geodesic equations were already derived by Monterde and Mu\~{n}oz Masqu\'e by a variational principle
(compare Thm. 7.1 and Cor. 7.1 in \cite{M-MM}). Apparently, the solutions of this equation were never studied in detail, excepting the observation that condition (i) above holds, already made in \cite{M-MM}. 
\vspace*{3mm}

Theorem \ref{uniqgeo} then gives existence and unicity of geodesics in a Riemannian supermanifold $\ca[N]$, with general
initial condition given by $\ca[J]$ in $\text{Mor}(\ca[S],T\can)$, with $\ca[S]$ being an arbitrary auxiliary supermanifold.
\vspace*{3mm}

We also recall there the relation between locally free sheaves and supervector bundles
as objects in the category of supermanifolds, mainly in order to have a precise grip on (co-)tangent bundles.
\vspace*{5mm}

In Section 3, we first remind the reader of the class of superfunctions on the total space of a supervector bundle
that are polynomial on the fibers, allowing us to define then a natural energy function $H$ on the cotangent
bundle $T^*\ca$ of a Riemannian supermanifold $(\ca , g)$. We give intrinsically (and in coordinates) the canonical 
symplectic form on $T^*\ca$ and explicit the Hamiltonian vector field $X_H$ associated to $H$. We proceed then to 
show that ``integral supercurves" of this vector field (i.e., curves given by the geodesic flow) are in bijection with geodesics on $(\ca , g)$, see Theorems \ref{thetheorem} and \ref{rthetheorem}.
\vspace*{5mm}

The last section, 4, gives applications of the geodesic flow morphism constructed in the preceding section. 
First, we show the existence of an exponential map generalizing the ungraded case and such that 
$\exp_q: T_q\ca \to \ca$, for $q$ in the body of $\ca$, is a local diffeomorphism (Theorem \ref{expdiffeo}).
Finally, we prove the result that
an isometry of a connected Riemannian supermanifold fixing a point of the body and having the identity
as its tangent map in this point, is already the identity morphism of the supermanifold. This implies immediately the unicity
of the ``geodesic symmetries" in the sense of Riemannian symmetric spaces, fixing a point in the body and having 
there minus the identity as their tangent map.

\section{Riemannian metrics on supermanifolds}
In this section we briefly recall  several fundamental notions, giving ideas of proof of certain folkloristic facts (Prop. \ref{thchris}, \ref{vectoralongmap} and \ref{aptang}) only when we did not find a reference in the literature. We define supergeodesics (Def. \ref{defgeo})
differently from \cite{Goertsches}, but our definition also reduces to the standard definition in the ungraded case.
We end this section with the first important new result of the article, namely existence and unicity of a supergeodesic for a given initial condition (Thm. \ref{uniqgeo}).\\

We use throughout this article the sheaf-theoretic approach to supermanifolds going back to Kostant \cite{Kostant} and Leites \cite{Leites}, thus basing our work on the following:

\begin{de}A \textbf{supermanifold $\ca$ of dimension $m|n$} is a locally ringed space $(M,\om)$, where $\om$ is a sheaf of commutative $\bb[Z]_2$-graded algebras on $M$, such that
\begin{enumerate}[(i)]
 \item $M$ is a Hausdorff second countable topological space,
 \item every point $m\in M$ has a neighborhood $U$ such that the ringed space $(U,\om[M]_{|U})$ is isomorphic to a superdomain $\ca[V]^{m|n}=(V,\ca[O]_{\sbb[m|n]|V})$ with $V$ open in $\sbb[m]$.
 \end{enumerate}
\end{de}

\begin{rem}If $\ca[X]=(X,\om[X])$ and $\ca[Y]=(Y,\om[Y])$ are ringed spaces, a \textbf{morphism} $\Phi:\ca[X]\to\ca[Y]$ is a pair $(\widetilde{\Phi},\Phi^*)$, where $\widetilde{\Phi}:X\to Y$ is a continuous mapping and
 $\Phi^*:\om[Y]\to\widetilde{\Phi}_*\om[X]$ is a collection of homomorphisms of rings $\Phi^*(U):\om[Y](U)\to\om[X](\widetilde{\Phi}^{-1}(U))$, for each open subset $U$ of $Y$, compatible with the restriction mappings. If the ringed spaces are locally ringed spaces, i.e., all stalks have a unique maximal ideal, the morphism $\Phi^*_x:\ca[O]_{\ca[Y],\widetilde{\Phi}(x)}\to \ca[O]_{\ca[X],x}$ is supposed to be local for all $x\in X$, i.e., $\Phi^*_x(\mathfrak{m}_{\widetilde{\Phi}(x)})\subset\mathfrak{m}_x$.
\end{rem}
We remark that the datum of such a collection $\Phi^*$ of homomorphisms is equivalent to $\Phi^*:\widetilde{\Phi}^{-1}\om[Y]\to\om[X]$ a collection of homomorphisms of rings $\Phi^*(V):(\widetilde{\Phi}^{-1}\om[Y])(V)\to\om[X](V)$, for each open subset $V$ of $X$, compatible with the restriction mappings. Recall that the inverse image $\widetilde{\Phi}^{-1}\om[Y]$ of $\om[Y]$ is defined as 
\[
 (\widetilde{\Phi}^{-1}\om[Y])(V):=\displaystyle\lim_{\substack{ \longrightarrow \\ W\supset \widetilde{\Phi}(V)\\ W\text{ open in }Y }}\om[Y](W).
\]
This equivalence is due to the fact that  $\widetilde{\Phi}_*$ and $\widetilde{\Phi}^{-1}$ are adjoint functors, i.e.
\[
 \text{Hom}\big(\om[Y],\widetilde{\Phi}_*\om[X]\big)\overset{\cong}{\longrightarrow} \text{Hom}\big(\widetilde{\Phi}^{-1}\om[Y],\om[X]\big).
\]
\begin{prop}Given supermanifolds $\cat$, $\cant$ the following holds:
 \begin{enumerate}[(i)]
  \item there is a morphism of sheaves $\beta:\om\to\ca[C]^\infty_M$, called the reduction mapping, such that the following sequence:
\[
0\longrightarrow\ca[O]_{\ca}^{\text{nil}}\longrightarrow\om\overset{\beta}{\longrightarrow}\ca[C]^\infty_M\longrightarrow 0
\]
is exact, where $\ca[O]_{\ca}^{\text{nil}}$ is the sheaf of nilpotent elements of $\om$.
\item if $\Phi:\ca\to\can$ is a morphism of supermanifolds, then we have the following commutative diagram
\[
 \begin{diagram}
  \node{\on}\arrow{e,t}{\Phi^*}\arrow{s,l}{\beta}\node{\widetilde{\Phi}_*\om}\arrow{s,r}{\beta}\\
\node{\ca[C]^\infty_N}\arrow{e,t}{\left(\widetilde{\Phi}\right)^{*}}\node{\widetilde{\Phi}_*\ca[C]^{\infty}_M\:\:.}
 \end{diagram}
\]
 \end{enumerate}
\end{prop}

\vspace*{5mm}
For every supermanifold $\ca=(M,\mathcal{O}_{\ca})$, one defines $\mathcal{T}_{\ca}(U):=\text{Der}_{\bb}(\mathcal{O}_{\ca}(U))$ for each open set $U$ in $M$. Here $\text{Der}_{\bb}(\mathcal{O}_{\ca}(U))$ denotes the $\mathcal{O}_{\ca}(U)$-module of derivations of $\mathcal{O}_{\ca}(U)$, where a homogeneous derivation $\varphi$ is by definition an element of \\ End$_{\bb}(\mathcal{O}_{\ca}(U))_0\cup\:$End$_{\bb}(\mathcal{O}_{\ca}(U))_1$ such that 
$$\varphi(ab)=\varphi(a)b+(-1)^{|\varphi||a|}a\varphi(b)\text{ for all }a,b\in \mathcal{O}_{\ca}(U).$$
$\tm$ is called the \textbf{tangent sheaf} and has the following structure (compare, e.g. \cite{Constantinescu}, Satz 4.36 for he proof):

\begin{theo}
$\mathcal{T}_{\ca}$ is a locally free sheaf of $\mathcal{O}_{\ca}$-modules of dimension $(m,n)$.
\end{theo}
The sections of $\mathcal{T}_M$ are called \textbf{vector fields}. The $\bb$-vector space $\tm(U)$ has a natural Lie superalgebra structure given by the following \textbf{bracket} $[.,.]$
\[
 [X,Y]f:=X(Yf)-(-1)^{|X||Y|}Y(Xf),
\]
satisfying the following graded Jacobi identity for each triple $X,Y,Z$ of vector fields:
\[
 [X,[Y,Z]]=[[X,Y],Z]+(-1)^{|X||Y|}[Y,[X,Z]].
\]
On a supermanifold one also has the \textbf{numerical tangent space}, defined as follows:\\

For every point $p$ in $M$, the numerical tangent space $T_p^{\text{num}}\ca$ of $\ca$ at $p$ is the space spanned by the homogeneous derivations $\varphi:\om_{,p}\to \bb$ satisfying
\[
 \varphi(fg)=\varphi(f)\widetilde{g}(p)+(-1)^{|\varphi||f|}\widetilde{f}(p)\varphi(g),
\]

where $\displaystyle \om_{,p}:=\lim_{\overrightarrow{U\ni p}}\om(U)$ is the stalk of $\om$ at $p$ and $\widetilde{f}(p)=\beta(f)(p)$.
$\\ \\$
For any $p \in U$, there is a natural reduction map $\mathcal{T}_{\ca}(U)\to T_p^{\text{num}}\ca$ sending a vector field $X\in\tm(U)$ to its value $X^{\text{num}}(p)$ at $p$.
Remark that, contrary to the classic case, a vector field is not determined by the collection of its numerical values at all points. 
\\

There is a natural real vector bundle of rank $m+n$ over $M$ having the numerical tangent spaces as its fibers, canonically splitting as the direct sum of $(T^{\text{num}}\ca)_0\to M$ and $(T^{\text{num}}\ca)_1\to M$. Then $X^{\text{num}}$ is a smooth section of $T^{\text{num}}\ca\to M$. If $\ca =
(M,\Gamma^\infty_{\Lambda E^*})$ with $E$ a smooth real vector bundle over $M$, the ``odd part" of the numerical tangent bundle, $(T^{\text{num}}\ca)_1$, is in fact isomorphic to $E$.
\\

We use again the letter $\beta$ to denote the reduction map to the numerical bundle composed with the projection to the even part of the tangent space:
\[
 \beta:\tm\to \left(M,\Gamma^\infty_{(T^{\text{num}}\ca )_0}  \right),
\]
where $ \left(M,\Gamma^\infty_{(T^{\text{num}}\ca )_0}  \right) $ is naturally isomorphic to the sheaf $(M,\Gamma^\infty_{TM})$ of vector fields on $M$.\\
For $f\in\om(U)$ and $X\in\tm(U)$, we will often write $\widetilde{f}$ and $\widetilde{X}$ instead of $\beta(f)$ and $\beta(X)$, respectively.
\begin{de}
The \textbf{cotangent sheaf} of a supermanifold $\ca$ is defined as the dual $\Omega^1_{\ca}$ of $\tm$.
\end{de}
\begin{prop}
The $\om$-module sheaf $\Omega^1_{\ca}$ is locally free of rank $(m,n)$.
\end{prop}
\begin{de} A \textbf{graded Riemannian metric} on a supermanifold $\ca$ is a graded symmetric
even non-degenerate $\ca[O]_{\ca}$-linear morphism of sheaves
\[
g:\ca[T]_{\ca}\otimes_{\om}\ca[T]_{\ca}\to\ca[O]_{\ca},
\]
where non-degeneracy means that the mapping $X\mapsto g(X,.)$ is an isomorphism from $\ca[T]_{\ca}$ to $\Omega^1_{\ca}.$
\end{de}
A supermanifold equipped with a graded Riemannian metric is called a \textbf{Riemannian supermanifold}.\\

\textbf{Remarks.} Let $q\in M$, the morphism $g$ defines a scalar superproduct $g_q$ on the real supervector space $T^{\text{num}}_q\ca$ as follows: Let $u,v$ be two numerical tangent vectors in $T^{\text{num}}_q\ca$. Choose vector fields $X,Y\in\tm(U)$ where $U$ is an open neighborhood of $q$ in $M$, fulfilling $X^{\text{num}}(q)=u$ and $Y^{\text{num}}(q)=v$. We then set $g_q(u,v):=\widetilde{g( X,Y)} (q).$ In general, this family of scalar superproducts does not determine the graded Riemannian metric.
By restriction to the even numerical vector fields this scalar superproduct induces a pseudo-Riemannian metric $\widetilde{g}$ on $M$, whereas the restriction to the odd numerical vector fields gives a symplectic structure
on the real vector bundle $(T^{\text{num}}\ca )_1\to M$. It follows notably that the odd dimension of a Riemannian supermanifold is always even.
\\

Our goal is to define the notion of supergeodesics, as ``autoparallel curves" with respect to a connection. Let thus
$\ca=(M, \om )$ be a supermanifold and $\mathcal{E}$ a locally free sheaf of $\om$-modules of finite rank on M.
\begin{de}
A \textbf{connection} on $\mathcal{E}$ is an even morphism $\nabla:\mathcal{E}\to \Omega^1_{\ca}\otimes_{\om}\mathcal{E}$ of sheaves of $\bb$-supervector spaces that satisfies the Leibniz rule
\[
 \nabla(f\cdot v)=df\otimes v+f\cdot \nabla (v)
\]
for all sections $f\in\om(U)$ and $v\in\mathcal{E}(U)$ where $df(X):=\sig[|X||f|]X(f)$ for all $X\in\tm(U)$, and on all open subsets $U$ of $M$.
\end{de}
Given a vector field $X$ on $\ca$, we have a map $\nabla_X:\ca[E]\to\ca[E]$ defined via the canonical pairing $\langle.,.\rangle$ between the tangent and the cotangent sheaf, $\nabla_X (v)=\langle X,\nabla(v)\rangle$. Then $\nabla_X$ satisfies:
\[
  (\nabla_X)(f \cdot v) = X(f)\cdot v + (-1)^{|X||f|} f \cdot (\nabla_X) (v)\quad\text{and}\quad |(\nabla_X)( v)| = |X| + |v|.
\]
    In the case $\mathcal{E} = \ca[T]_M$, we speak of a \textbf{connection on} $\ca$.\\
\begin{de}
We define
the \textbf{torsion} of a connection $\nabla$ on $\ca$ by
\[
 T_\nabla (X, Y ):= \nabla_X Y-(-1)^{|X||Y |}\nabla_Y X - [X, Y ].
\]
\end{de}
\begin{de}
    If $\ca$ is equipped with a graded Riemannian metric $g$, a connection $\nabla$ is called \textbf{metric} if the connection respects the metric in the following sense
\[
 \nabla g=0\text{, i.e. }  X(g(Y\otimes Z)) = g(\nabla_X Y\otimes Z) + (-1)^{|X||Y |}g( Y\otimes \nabla_X Z)\quad\forall X,Y,Z\in\tm(U) . 
\]
\end{de}
We also have the natural graded analogue of \textbf{Christoffel symbols}: If $(q_i )$ is a system of
coordinates on $U\subset M$, the expansion
\[
 \nabla_{\pai}\paj=\sum_k\Gamma_{ij}^k\pak
\]
gives elements $\Gamma_{ij}^k\in \ca[O]_M (U )$ of parity $|\Gamma_{ij}^k|=|q_i|+|q_j|+|q_k|$.

\begin{theo}\label{formula-levi}
  On a supermanifold $\ca$ with a graded Riemannian metric $g$, there
exists a unique torsion free and metric connection $\nabla$, which will be called the
\textbf{Levi-Civita connection} of the metric. \end{theo}
The proof follows immediately upon interpreting the usual six-term-formula in a graded sense and sheaf-theoretically.\\

A connection on $\ca$ induces a connection on the vector bundle $TM\to M$ by reduction. Indeed, there exists a unique connection $\widetilde{\nabla}$ on $M$ such that 
\[
 \widetilde{\nabla}_{\beta(X)}(\beta(Y))=\beta(\nabla_X Y)\quad\forall X,Y\in \ca[T]_M(U)\text{ and for all }U\text{ open of }M.
\]

In our case, the Levi-Civita connection associated to a graded Riemannian metric induces the usual Levi-Civita
connection of the induced metric $\widetilde{g}$ on $M$. Moreover, in coordinates $(q_i)$, the reduction $\beta(\Gamma^k_{ij})$ of a
Christoffel symbol $\Gamma^k_{ij}$ is the Christoffel symbol of the Levi-Civita connection on $M$ with respect to the corresponding coordinates on $M$.
\\

In the classical case, we have an explicit formula for Christoffel symbols. Adding appropriate signs this formula holds true in the super context:
\begin{prop}\label{thchris}In coordinates $(q_i)$ on $U$,
 the Christoffel symbols $\Gamma^k_{ij}$ of the Levi-Civita connection on a Riemannian supermanifold are given by
\begin{equation}
  \displaystyle\Gamma^k_{ij}=\frac{1}{2}\sum_l\left[ \pa[q_i]g_{jl}+\sig[|q_i||q_j|]\pa[q_j]g_{il}-\sig[|q_l|(|q_i|+|q_j|)]\pa[q_l]g_{ij} \right]g^{lk},\label{chritoffel}
\end{equation}
where $g_{ij}:=g(\pa[q_i],\pa[q_j])\in\om(U)$ and $g^{lk}\in\om(U)$ are defined by $\displaystyle\sum_k g^{ik}g_{kj}=\delta_{i,j}$.
\end{prop}

\textbf{Proof.} Direct computation using the six-term-formula that defines the Levi-Civita connection:
\[
 2 g( \nabla_X Y\otimes Z) = Xg( Y\otimes Z) - (-1)^{|Z|(|X|+|Y |)} Zg( X\otimes Y)
             + (-1)^{|X|(|Y |+|Z|)}Y g( Z\otimes X ) +\]
\[ g( [X, Y ]\otimes Z)
             - (-1)^{|X|(|Y |+|Z|)}g( [Y, Z]\otimes X)
             + (-1)^{|Z|(|X|+|Y |) }g( [Z, X]\otimes Y) 
\]for all $X,Y$ and $Z$ homogeneous vector fields.
\hfill $\square$\\

A connection on $\ca$ allows to define a covariant derivative of vector fields along a curve. We develop
these notions here in some detail for supermanifolds.
\begin{de}
Let $\Phi=(\widetilde{\Phi},\Phi^*):\ca\to \ca[N]$ be a morphism of supermanifolds and $U$ an open set of $N$. A \textbf{vector field along} $\Phi:\ca\to\ca[N]$ on $U$ is a morphism of supervector spaces 
\[
 X:\ca[O]_{\ca[N]}(U)\to (\widetilde{\Phi}_*\ca[O]_{\ca})(U)
\]
such that its homogeneous components $X_0$, $X_1$  satisfy the following derivation property
\[
 X_\alpha(fg)=X_\alpha (f)\Phi^*(g)+(-1)^{\alpha \cdot|f|}\Phi^*(f)X_\alpha(g)\quad\text{for }\alpha=0,1\text{ and }f,g\in\ca[O]_{\ca[N]}(U).
\]
\end{de}
 The set $\text{Der}_{\bb}\left(\om[N](U),(\widetilde{\Phi}_*\om)(U)\right)$ of all vector fields along $\Phi$ will also be denoted by $\ca[T]^{\Phi}_N(U)$. Remark that $\tphi_N$ is naturally a sheaf of $\widetilde{\Phi}_*\om$-modules over $N$. \\

Recall that a morphism of supermanifolds
can be viewed as a morphism $\Phi^*:\on\to\widetilde{\Phi}_*\om$ of sheaves over $N$ or as a morphism $\Phi^*:\widetilde{\Phi}^{-1}\on\to\om$ of sheaves over $M$. Accordingly, vector fields along $\Phi$ can also be viewed as a sheaf of $\om$-modules over $M$: $\tphiup_M:=\text{Der}_{\bb}\left(\widetilde{\Phi}^{-1}\on,\om\right)$.\\

There are two standard ways of constructing vector fields along $\Phi$. If $X$ is a vector field on $\ca[N]$, then $$ \widehat{X}:=\Phi^*\circ X$$ is a vector field along $\Phi$, and if $Y$ is a vector field on $\ca$, we get the following vector field along $\Phi$: $$ (T\Phi^*)(Y):=Y\circ \Phi^*.$$
We thus obtain two morphisms of sheaves over $M$, namely $\widetilde{\Phi}^{-1} \tm[N]\overset{\Phi^*}{\longrightarrow}\ca[T]_M^{\Phi}$ and $\tm[M]\overset{T\Phi^*}{\longrightarrow}\ca[T]_M^\Phi$, and their analogues over $N$.
Since knowledge of this sheaves will be useful in the sequel, we give a more explicit description.
\begin{prop}\label{vectoralongmap}$\\$
\vspace{-5mm}
\begin{enumerate}[(i)]
\item  $\ca[T]_N^\Phi$ is a locally free sheaf of $\widetilde{\Phi}_{*}\om$-modules over $N$ of the same rank as $\tm[N]$. More precisely, we have an isomorphism
\[
 \ca[T]_N^{\Phi}\overset{\cong}{\longrightarrow}\widetilde{\Phi}_*\om\otimes_{\om[N]}\tm[N].
\]
\item  $\ca[T]^\Phi_M$ is a locally free sheaf of $\om$-modules over $M$ of the same rank as $\tm[N]$. More precisely, we have an isomorphism
\[
 \ca[T]^{\Phi}_M\overset{\cong}{\longrightarrow}\om\otimes_{\widetilde{\Phi}^{-1}\om[N]}\widetilde{\Phi}^{-1}\tm[N].
\]
\end{enumerate}
\end{prop}
\textbf{Proof.} We will prove only \textit{(i)}. Part \textit{(ii)} then follows from the fact that $\widetilde{\Phi}_* $ and $\widetilde{\Phi}^{-1}$ are adjoint functors.\\
We will construct two isomorphisms $Q$ and $S$ as below:
\[
\begin{diagram}
\node{\text{Der}_{\bb}\left(\om[N],\widetilde{\Phi}_*\om\right) }\arrow{e,tb}{Q}{\cong}\node{\text{Hom}_{\om[N]}\left(\Omega_{\ca[N]}^1,\widetilde{\Phi}_*\om\right)}\arrow{e,tb}{S}{\cong}\node{\widetilde{\Phi}_*\om\otimes_{\om[N]}\tm[N].}
\end{diagram}
\]
\textit{Step 1}:\\
We define $Q$ on $U$, a chart domain of $\ca[N]$ with local coordinates $(q_i)$.
Let $D$ be a vector field along $\Phi$ on $U$, $Q(U)(D)$ is then the morphism of $\om[N](U)$-module such that
\[
 Q(U)(D)(dq_i)=\sig[|i|]D(q_i)\quad\forall i,
\]
where $dq_i(\pa[q_j]):=\delta_{i,j}$.\\

Thus if $\omega=\displaystyle\sum_i\omega_i\cdot dq_i$, we have $  Q(U)(D)(\omega)=\displaystyle\sum_i\sig[|D||\omega_i|+|i|]\Phi^*(\omega_i)\cdot D(q_i)$.
It follows that $Q(U)$ is a morphism of $(\widetilde{\Phi}_*\om)(U)$-modules and $Q$ is a $\widetilde{\Phi}_*\om$-module morphism.

\textit{Step 2}:\\
Similarly, we define $S$ on a chart domaine $U$. Let $\varphi\in\text{Hom}_{\om[N](U)}\left(\Omega_{\ca[N]}^1(U),\widetilde{\Phi}_*\om(U)\right)$, then
\[
S(U)(\varphi)=\displaystyle\sum_i\sig[|i|] \varphi(dq_i)\otimes_{\om[N](U)}\pa[q_i],\quad\text{where}\quad\pa[q_i]:=\displaystyle\frac{\partial}{\pa[q_i]}.
\]
A direct computation shows that both maps $Q$ and $S$ are isomorphisms of $\widetilde{\Phi}_*\om$-modules.\hfill $\square$\\

\textbf{Remark.} The proof of Proposition \ref{vectoralongmap} shows that the isomorphism $L:=S\circ Q$ looks as follows in local coordinates on $U$, for a $D\in\tphi_N(U)$:
\[
\begin{array}{ccl}
 L(D)     &=&\displaystyle\sum_i\sig[|i|] Q(D)(dq_i)\otimes\pa\\
     &=&\displaystyle\sum_i D(q_i)\otimes\pa[q_i].
\end{array}
\]
Furthermore, if $(q_i)$ are coordinates on $\can$, then $(\hat{\partial}_{q_i})$ is a basis of $\ca[T]_N^\Phi(U)$ and we have
\[
 T\Phi^*(Y)=Y\circ\Phi^*=\sum_i Y(\Phi^*q_i)\hat{\partial}_{q_i}
\]
for all vector fields $Y$ on $M$.\\

Given a connection on $\ca[N]$ and a morphism $\Phi:\ca\to\ca[N]$, we introduce the associated pull-back connection in order to get a covariant derivative along a curve. One has, in fact, the statement below:
\begin{propdef}\label{pullbackcon} There exists a unique connection $\widehat{\nabla}$ on $\tphi_N$ such that the diagram
\[
\begin{diagram}
\node{\tn}\arrow{e,t}{\nabla}\arrow{s,l}{\Phi^*}\node{\Omega_{\ca[N]}^1\otimes_{\om[N]}\tn}\arrow{s,r}{\Phi^*}\\
\node{\tphi_N}\arrow{e,t}{\widehat{\nabla}}\node{\widetilde{\Phi}_*\Omega_{\ca}^1\otimes_{\widetilde{\Phi}_*\om}\tphi_N}
\end{diagram}
\]
commutes. This connection is usually called the \textbf{pullback of $\nabla$ via} $\Phi$ and is denoted by $\Phi^*\nabla$.
\end{propdef}
\textbf{Remark. }The vertical arrow is define as follows.
First there is a natural map from $\Omega_{\ca[N]}^1$ to $\widetilde{\Phi}_*\Omega_{\ca}^1$ induced by $\Phi^*$ which is also denoted by $\Phi^*$.
Given $U$ an open set in $N$ and $\omega\in\Omega_{\ca[N]}^1(U)$, we get an element $\upsilon$ of $\widetilde{\Phi}_*\Omega_{\ca}^1(U)$ :
\[
 \upsilon(X):=\widehat{\omega}(T\Phi^*(X))\text{ for }X\in(\widetilde{\Phi}_*\tm)(U),
\]
where $\widehat{\omega}$ is defined by
\[
 \begin{array}{cccc}
  \widehat{\omega}:&\widetilde{\Phi}_*\om(U)\otimes_{\on(U)}\tn(U)&\longrightarrow&\widetilde{\Phi}_*\om(U)\\
                   &f\otimes Y                                    &\longmapsto    &\sig[|f||\omega|]f\cdot\Phi^*(\omega(Y))
 \end{array}
\]
and the isomorphism $\ca[T]_N^{\Phi} \cong\widetilde{\Phi}_*\om\otimes_{\om[N]}\tm[N]$.
Using $\Phi^*:\on\to\widetilde{\Phi}_*\om$ the map $\Phi^*:\Omega_{\ca[N]}^1\longrightarrow\widetilde{\Phi}_*\Omega_{\ca}^1$ turns out to be a morphism of $\on$-modules, that can be tensored with $\Phi^*:\tn\to\tphi_N$ to obtain the arrow
\[
 \Phi^*:\Omega_{\ca[N]}^1\otimes_{\om[N]}\tn\longrightarrow\widetilde{\Phi}_*\Omega_{\ca}^1\otimes_{\widetilde{\Phi}_*\om}\tphi_N.
\]
\textbf{Proof of Proposition \ref{pullbackcon}.} The proof follows immediately upon interpreting the usual proof in the classical case of ungraded manifolds sheaf-theoretically.\hfill$\square$\\

\textbf{Remark.} The diagram in Proposition/Definition \ref{pullbackcon} can be equivalently replaced by the following diagram:
\[
\begin{diagram}
\node{\widetilde{\Phi}^{-1}\tn}\arrow{e,t}{\nabla}\arrow{s,l}{\Phi^*}\node{\widetilde{\Phi}^{-1}\Omega_{\ca[N]}^1\otimes_{\widetilde{\Phi}^{-1}\om[N]}\widetilde{\Phi}^{-1}\tn}\arrow{s,r}{\Phi^*}\\
\node{\tphi_M}\arrow{e,t}{\widehat{\nabla}}\node{\Omega_{\ca}^1\otimes_{\om}\tphi_M}
\end{diagram}
\]
We can now define the covariant derivative along a supercurve $\Phi:\sbb\to\ca[N]$. Let $(\ca[N],g)$ be a Riemannian supermanifold and $\nabla$ its Levi-Civita connection. Consider $\Phi^*\nabla$, the pullback connection of $\nabla$ under $\Phi$. 
\begin{de}\label{defcov}The \textbf{even covariant derivative along} $\Phi:\sbb\to\ca[N]$ is defined as
\[
 \displaystyle\frac{\nabla}{dt}:=(\Phi^*\nabla)_{\pa[t]}:\tphi_N\longrightarrow\tphi_N,
\]
induced by the duality pairing of $\pa[t]$ with elements of $\Omega_{\sbb}^1$ where $(t,\theta)$ are the cananical coordinates on $\sbb\cong\bb\times\sbb[0|1]$.\\
\end{de}

Note that if $(q_i)$ are local coordinates on $U\subset N$, we have for all vector fields $X\in\tphi_N(U)$
\[
 \begin{array}{ccl}
 \displaystyle\frac{\nabla}{dt}(X)&=& \displaystyle \sum_k\left[\pa[t]X(q_k)  + \sum_{i,j}  X(q_i)\pa[t]\Phi^*(q_j)\Phi^*(\Gamma_{ji}^k)\right] \otimes\pa[k].
 \end{array}
\]
Analogously, we can define the covariant derivative along a curve with respect to the odd coordinate $\theta$ on $\sbb$:
\[
  \displaystyle\frac{\nabla}{d\theta}:=(\Phi^*\nabla)_{\pa[\theta]}:\tphi_N\longrightarrow\tphi_N,
\]
induced by the duality pairing of $\pa[\theta]$ with elements of $\Omega_{\sbb}^1$.\\
In local coordinates, we have a similar formula with an additional sign:
\[
 \begin{array}{ccl}
 \displaystyle\frac{\nabla}{d\theta}(X)&=&\displaystyle \sum_k\left[\pa[\theta]X(q_k)  + \sum_{i,j}  \sig[|X|+|q_i|]X(q_i)\pa[\theta]\Phi^*(q_j)\Phi^*(\Gamma_{ji}^k)\right] \otimes\pa[k]
 \end{array}
\]

We can now give a natural definition of  supergeodesics.
\begin{de}\label{defgeo}A \textbf{supergeodesic} on $(\ca[N],\nabla)$, a supermanifold with a connection $\nabla$, is a supercurve $\Phi:\sbb\cong\bb\times\sbb[0|1]\to\ca[N]$ such that
\[
 \displaystyle\frac{\nabla}{dt}\big(T\Phi^*(\pa[t])\big)=0.
\]
\end{de}
The preceding condition is locally equivalent to
\begin{eqnarray}
  \pa[t]^2 \Phi^*(q_k)  + \sum_{i,j}  \pa[t]\Phi^*(q_i)\pa[t]\Phi^*(q_j)\Phi^*(\Gamma_{ji}^k)=0\quad\forall k, \label{equageodesic}
\end{eqnarray}
where $(q_k)$ are local coordinates on $\ca[N]$.
Note that the equation (\ref{equageodesic}) restricted to those $k$ corresponding to even coordinates are the usual geodesic equations of the underlying classical manifold equipped with the reduced connection $\widetilde{\nabla}$ on $M$.\\

\textbf{Remarks.} (1) Of course, in general a supergeodesic is only defined on an open subsupermanifold of $\sbb$ with body an open interval in $\bb$. Slightly abusing, we will often write $\sbb$ meaning such an open subsupermanifold of it.\\
(2) One can easily generalize Definitions \ref{defcov} and \ref{defgeo} to the case $\bb\times\ca[S]$, $\ca[S]$ an arbitrary supermanifolds replacing $\sbb[0|1]$.\\

A different definition was given by Goertsches (\cite{Goertsches}):
a \textbf{supergeodesic in the sense of Goertsches} is a supercurve $\Phi:\sbb\to\ca[N]$ such that
\begin{eqnarray}
 \left(\frac{\nabla}{dt}T\Phi^*(\pa[t])\right)^{\text{num}}= \left(\frac{\nabla}{dt}T\Phi^*(\pa[\theta])\right)^{\text{num}}=0. \label{equageodesicgeortches}
\end{eqnarray}
Remark that such a curve automatically satisfies the following equations:
\[ 
\left(\frac{\nabla}{d\theta}T\Phi^*(\pa[t])\right)^{\text{num}}= \left(\frac{\nabla}{d\theta}T\Phi^*(\pa[\theta])\right)^{\text{num}}=0.
\]
Locally, (\ref{equageodesicgeortches}) is equivalent to 
\[
  f_k''+\sum_{\substack{i\text{ even,}\\j\text{ even}}}f'_i\cdot f'_j\cdot \Phi^*(\Gamma_{ji}^k)=0\quad\forall k\text{ even},
\]
and
\[
f'_\delta+\sum_{\substack{i\text{ even,}\\ \beta\text{ odd}}}f_\beta \cdot f'_i\cdot \Phi^*(\Gamma^\delta_{i\beta})=0\quad\forall \delta\text{ odd},
\]
where $f_i,f_\delta\in\ca[C]_{\bb}^\infty(\widetilde{\Phi}^{-1}(U))$ are functions defined as follows:
\[
 \Phi^*(q_i)=f_i\quad\forall i\text{ even and }\Phi^*(q_\beta)=f_\beta\cdot \theta\quad\forall \beta\text{ odd.}
\]
We call here and in the sequel an index $(i,\delta,k,\beta,\ldots)$ ``even" resp. ``odd" if the corresponding coordinate has the property.\\

In particular, if $\ca[N]=\sbb[1|2]$ is equipped with a flat metric (i.e., $\Gamma_{ij}^{k}=0$), (\ref{equageodesicgeortches}) is equivalent to
\[
  f_k''=0\quad\forall k\text{ even, and }f'_\delta=0\quad\forall \delta\text{ odd}.
\]
On the other hand, the geodesic equation (\ref{equageodesic}) is equivalent to
\[
  f_k''=0\quad\forall k\text{ even, and }f''_\delta=0\quad\forall \delta\text{ odd}.
\]
The equations (\ref{equageodesic}) and (\ref{equageodesicgeortches}) are thus obviously non-equivalent, moreover initial conditions for these two notions of geodesics are not the same.\\

Let us explain in some detail what an initial condition is in the super context. For this we need to recall supervector bundles and the tangent map (in the category of supervector bundles) coming from a morphism of supermanifolds.\\

If $\ca[X]=(X,\om[X])$ is a supermanifold and $V$ an open subset of $X$, we denote by $\ca[X]_{|V}$ the open subsupermanifold given by $(V,\ca[O]_{\ca[X]|V})$.

\begin{de} Let $\ca[V]$ and $\ca$ be supermanifolds and let $\pi:\ca[V]\to\ca$ be a surjective submersion.\\
 $\ca[V]$ is called a \textbf{supervector bundle} with base $\ca$ and typical fiber $\sbb[p|q]$ if
\begin{enumerate}[(i)]
 \item  for all $x\in M$, there exists a local trivialization $(U,\varphi)$, i.e. there is an open neighborhood $U$ of $x$ and an isomorphism $\varphi:\ca[V]_{|\widetilde{\pi}^{-1}(U)}\to\sbb[p|q]\times\ca_{|U}$ such that the following diagram commutes
\[
 \begin{diagram}
  \node{\ca[V]_{|\widetilde{\pi}^{-1}(U)}}\arrow{e,t}{\varphi}\arrow{s,l}{\pi}\node{\sbbpq\times\ca_{|U}}\arrow{sw,r}{\text{proj}_{1}}\\
\node{\ca_{|U}}
 \end{diagram}
\]
\item and there exists a covering of $M$ by local trivializations $\{(U_\alpha,\varphi_\alpha)\}$ as in (i) such that for all $\alpha,\beta$ with $U_{\ab}:=U\al\cap U\be$ non empty set, the transition function 
$\varphi_{\alpha\beta}:=\varphi\be\circ(\varphi\al)^{-1}:\sbb[p|q]\times\ca_{|U_{\ab}}\to\sbb[p|q]\times\ca_{|U_{\ab}}$ 
is determined by an element $\phi_{\alpha\beta}\in GL_{(p,q)}(\om(U_{\ab}))$, that is 
\[
\begin{array}{ccl}
 (\varphi_{\ab})^*(f)&=&f\quad\forall f\in \om(U_{\ab}),\\
  (\varphi_{\ab})^*(v_i)&=&\displaystyle\sum_j v_j \cdot(\phi_{\alpha\beta})_{ji}          \quad\forall i,
 \end{array}
\]
where $(v_i)$ are the canonical coordinates on $\sbbpq$.
\end{enumerate}
\end{de}
\begin{theo}\label{modtosvec}(\cite{Schmitt}) Let $\ca$ be a supermanifold. The category of locally free sheaves of $\om$-modules of finite rank is equivalent to the category of supervector bundles over $\ca$.
\end{theo}
Now let $T\ca=(TM,\ca[O]_{T\ca})$ and $T^*\ca=(T^*M,\ca[O]_{T^*\ca})$ be the supervector bundles associated to the locally free sheaves $\tm$ and $\Omega_{\ca}^1$ respectively.
\begin{prop}\label{aptang}Each morphism of supermanifolds $\Phi:\ca\to\can$ induces a morphism of supermanifolds $T\Phi:T\ca\to T\can$ called \textbf{the tangent map} of $\Phi$.
\end{prop}
\textbf{Sketch of the proof.} Let $(q^{\ca}_i)$ and $(q^{\can}_j)$ be coordinates on $\ca$ and $\can$. They induce natural coordinates $(q^{\ca}_i,v^{\ca}_i)$ and $(q^{\can}_j,v^{\can}_j)$ which trivialize $T\ca$ and $T\can$, respectively. In these coordinates, the tangent map of $\Phi$ is given as follows
\begin{equation}
\begin{array}{ccl}
 (T\Phi)^*(q^{\can}_j)&=&\Phi^*(q^{\can}_j)\quad\forall j,\\
  (T\Phi)^*(v^{\can}_j)&=&\displaystyle\sum_i v^{\ca}_i\cdot\pa[{q^{\ca}_i}]\Phi^*(q^{\can}_j)          \quad\forall j.
 \end{array}\tag*{$\square$}
\end{equation}
Let $\Phi:\bb\times\ca\to\ca[N]$ be a morphism of supermanifolds. Using the tangent map $T\Phi:T\bb\times T\ca\to T\can$, the section $s_1:t\mapsto(1,t)$ of $T\bb$ and the canonical zero section $s_0:\ca\to T\ca$ of $T\ca$, we define $\pa[t]\Phi:\bb\times\ca\to T\can$, the \textbf{partial derivative of $\Phi$ in the direction of $\bb$}, as follows:
$$\pa[t]\Phi:=T\Phi\circ (s_1\times s_0):\bb\times \ca\to T\can.$$
In local coordinates, this amounts to
$$
\begin{array}{lcl}
(\pa[t]\Phi)^*(q_i^{\can}) & = & \Phi^*(q_i^{\can})\\
(\pa[t]\Phi)^*(v_i^{\can}) & = & \pa[t]\Phi^*(q_i^{\can}).
\end{array}
$$
We define for $\ca[X]$, $\ca[Y]$ supermanifolds and $x_0\in X$, $\text{inj}_{x=x_0}:\ca[Y]\to\ca[X]\times\ca[Y]$ by $\widetilde{\text{inj}_{x=x_0}}(y)=(x_0,y)$ and $(\text{inj}_{x=x_0})^*(f\otimes g)=\widetilde{f}(x_0)\cdot g$ for $f\in\om[X]$ and $g\in\om[Y]$. The pullback operation $(\text{inj}_{x=x_0})^*$ is often denoted by $\text{ev}_{x=x_0}$.

\begin{theo}\label{uniqgeo}
Let $\ca[J]\in\text{Mor}(\sbb[0|1],T\can)$. Then there exists a unique geodesic $\Phi:\sbb\cong\bb\times\sbb[0|1]\to\can$ such that $\pa[t]\Phi\circ\text{inj}_{t=0}=\ca[J]$.
\end{theo}
\textbf{Proof.} The ``initial condition" $\ca[J]$ renders the solution of the system of ODEs (\ref{equageodesic}) unique, since $\pa[t]\Phi\circ\text{inj}_{t=0}=\ca[J]$ reads in the local coordinates discussed before the theorem as follows: 
\begin{equation}
\text{inj}_{t=0}^*\circ\Phi^*(q_i)=\Phi^*(q_i)_{|t=0}=\ca[J]^*(q_i)\text{ and }\pa[t]\Phi^*(q_i)_{|t=0}=\ca[J]^*(v_i).\tag*{$\square$}
\end{equation}
\textbf{Remark. }The preceding theorem can be easily generalized to cover the situation discussed in the remark after Definition \ref{defgeo}: given a supermanifold $\ca[S]$ and a morphism $\ca[J]\in\text{Mor}(\ca[S],T\can)$, there exists a unique geodesic $\Phi:\bb\times\ca[S]\to\can$ such that $\pa[t]\Phi\circ\text{inj}_{t=0}=\ca[J]$.

\section{The geodesic flow}
In this section we first clarify the construction of superfunctions on supervector bundles coming from a locally free sheaf, as well as of the canonical two-form on the cotangent bundle of a supermanifolds (both already sketched in \cite{Kostant}).
We get a natural ``energy function" $H$ on the (co-)tangent bundle of a supermanifold and an associated Hamiltonian
vector field $X_H$. The flow and integral curves of the latter even vector field exist by results of Monterde and
Montesinos resp. Dumitrescu. We then prove that integral curves of $X_H$ on $T^*\ca$ are in natural bijection with geodesics on $\ca$, as asked for in the introduction. 

\begin{prop}\label{naturalmap}
Let $\ca=(M,\om)$ be a supermanifold and $\ca[F]$ be a locally free sheaf of $\om$-modules and let $\ca[V]=(V,\ca[O]_{\ca[V]})\overset{\pi}{\longrightarrow}\ca$ be the supervector bundle associated to $\ca[F]$. 
Then there exists a natural sheaf morphism $\widetilde{\pi}^{-1}\ca[F]^{*}\longrightarrow\ca[O]_{\ca[V]}\text{ over }V$ or equivalently a sheaf morphism $\ca[F]^{*}\longrightarrow\widetilde{\pi}_*\ca[O]_{\ca[V]}\text{ over }M.$
\end{prop}
\textbf{Proof.} The proof follows directly from the construction of the supervector bundle $\ca[V]$ associated to $\ca[F]$. \hfill $\square$\\
\textbf{Remark. }Section of $\widetilde{\pi}^{-1}\ca[F]^*$ over on open subset in $V$ are of course local superfunctions of $\ca[O]_{\ca[V]}(V)$ that are ``linear in the fiber directions".
\begin{cor} Let $\pi:\ca[V]\to\ca$ be the canonical projection. Then we have a natural morphism of sheaves
\[
\begin{diagram}
\node{\widetilde{\pi}^{-1}\left(\text{Sym}_{\om}\left(\ca[F]^*\right)\right)}\arrow[2]{e,t}{\chi}\arrow{se}\node[2]{\ca[O]_{\ca[V]}}\arrow{sw}\\
\node[2]{V}
\end{diagram}
\]
\end{cor}
\textbf{Proof.} By the preceding proposition we have a morphism of $\widetilde{\pi}^{-1}\om$-module sheaves $\widetilde{\pi}^{-1}\ca[F]^*\longrightarrow\om[V]$. Since $\om[V]$ is a sheaf of unital supercommutative superalgebras, we have a canonical extension to $\widetilde{\pi}^{-1}\left(\text{Sym}_{\om}(\ca[F]^*)\right)\cong\text{Sym}_{\widetilde{\pi}^{-1}\om}(\widetilde{\pi}^{-1}\ca[F]^*)\longrightarrow\om[V]$, denoted by $\chi$.\hfill $\square$ \\

In the case that $\ca[F]=\Omega^{1}_{\ca}$, we have, after identifying $(\Omega^1_{\ca})^*$ with $\tm$, notably the following morphism of $\om(M)$-modules:
\[
\begin{diagram}
 \node{\left(\text{Sym}_{\om}^2\left(\tm\right)\right)(M)}\arrow{e,t}{\chi(T^*M)}\node{\ca[O]_{T^*\ca}(T^*M).}
\end{diagram}
 \]
If $\ca$ is equipped with a graded Riemannian metric $g:\tm\otimes_{\om}\tm\to\om$, we have the isomorphisms $g^{\flat}:\tm\overset{\cong}{\longrightarrow}\Omega_{\ca}^1$ and $g^{\natural}:\Omega_{\ca}^1\overset{\cong}{\longrightarrow}\tm$ induced by the metric, where 
\[
 g^{\flat}(X)(Y):=g(X,Y)\quad\forall X,Y\in\tm(U).
\]
Using the above identification and the metric we have
\[
 \begin{diagram}
 \node{\Omega^1_{\ca}(M)\otimes\Omega^1_{\ca}(M)}\arrow{e,tb}{g^\natural\otimes g^\natural}{\cong}\node{\ca[T]_{\ca}(M)\otimes\ca[T]_{\ca}(M)}\arrow{e,t}{g}\node{\om(M),}
\end{diagram}
\]
which determines an element $h$ of $\left(\Omega_{\ca}^1(M)\otimes\Omega_{\ca}^1(M)\right)^*$. The standard theory of linear superalgebra gives the following canonical isomorphisms
$$
\left(\Omega_{\ca}^1(M)\otimes\Omega_{\ca}^1(M)\right)^*\cong\left(\Omega_{\ca}^1(M)\right)^*\otimes\left(\Omega_{\ca}^1(M)\right)^*\cong\tm(M)\otimes\tm(M).
$$
Let us now define $H:=\frac{1}{2}\chi(T^*M)(h)\in\mathcal{O}_{T^*\ca}(T^*M)$, where $h$ is viewed as a symmetric element of $\tm(M)\otimes\tm(M)$.
In local coordinates, we obtain 
\begin{eqnarray*}
  H &=& \displaystyle\frac{1}{2}\sum_{i,j}p_i\cdot g^{ij}\cdot p_j,
\end{eqnarray*}
where $(q_i,p_i)$ are coordinates on $T^*\ca$ induced by coordinates $(q_i)$ on $\ca$.\\
Recall that in the super context we have a natural notion of graded differential $k$-form $\Omega_{\ca}^{k}$ defined, e.g., in \cite{Kostant}.
 Moreover, we have a canonical exterior differentiation $d$ of bidegree $(1,0)$ of the $(\bb[Z]\times\bb[Z]_{2})$-graded commutative algebra $\Omega_{\ca}$ of graded differential forms.

\begin{de}A \textbf{symplectic supermanifold} is a supermanifold $\ca=(M,\om)$ together with a closed even non-degenerated two-form $\omega\in\Omega_{\ca}^2(M)_{0}$.
\end{de}
\begin{prop}Let $\ca$ be a supermanifold. Then the cotangent bundle $T^{*}\ca$ of $\ca$ is equipped with a canonical symplectic form $\omega\in\Omega_{T^{*}\ca}^{2}(T^{*}M)_{0}$. As in the classical case of ungraded manifold, $\omega=-d\alpha$ with a canonical potential $\alpha\in\Omega_{T^{*}\ca}^{1}(T^{*}M)_{0}$.
\end{prop}
\textbf{Proof.} Using the notations of Proposition \ref{vectoralongmap} and the remark after this proposition, we can define the canonical even one-form $\alpha$ on $T^{*}\ca$ as follows:
\[
 \begin{diagram}
  \node{\text{Der}_{\bb}(\ca[O]_{T^*\ca}(T^{*}M))}\arrow{e,t}{T(\pi^{T^*\ca})^*}\arrow{see,b}{\alpha} \node{\text{Der}_{\bb}(\om(M),\ca[O]_{T^*\ca}(T^{*}M))}\arrow{e,tb}{L}{\cong}\node{\ca[O]_{T^*\ca}(T^{*}M)\otimes_{\om(M)}\ca[T]_{\ca}(M)}\arrow{s,r}{1\otimes\chi}\\
\node[3]{\ca[O]_{T^*\can}(T^{*}M).}
 \end{diagram}
\]
Here $T^*\ca\overset{\pi^{T^*\ca}}{\longrightarrow}\ca$ is the canonical projection, $L$ is the map defined in the Proposition \ref{vectoralongmap} and $\chi$ is the morphism given by the Theorem \ref{naturalmap}. The projection $\pi^{T^*\ca}$ will simply be denoted by $\pi$ in the rest of this proof.\\
We obtain an element $\alpha\in\Omega_{T^{*}\ca}^{1}(T^{*}M)_{0}$ described as follows in the local coordinates $(q_i,p_i)$:
\begin{eqnarray*}
\alpha(X) & = & (1\otimes\chi)\circ L\circ T\pi^*(X)\\
	  & = & (1\otimes\chi)\circ L( X\circ\pi^*)\\
	  & = & (1\otimes\chi)\left(\displaystyle\sum_{i}  X\circ\pi^*(q_i)\otimes\pa \right)\\
	  & = & \displaystyle\sum_{i}  X( q_i)\cdot\chi(\pa),\text{ by identifying }\pi^*(q_i)\text{ with }q_i,\\
	  & = & \displaystyle\sum_{i}  X(q_i)\cdot p_i\\
	  & = & \displaystyle\sum_{i}\sig[|p_i|(|q_i|+|X|)]p_i\cdot  X( q_i)\\
	  & = & \displaystyle\sum_{i}\sig[|p_i|]p_i\cdot d(q_i )(X)\quad\forall X\in\ca[T]_{T^{*}\ca}(T^*M).
\end{eqnarray*}
We thus locally have $\alpha=\displaystyle\sum_{q_i}\sig[|q_i|]p_i\cdot d( q_i )$.\\

Now we set $\omega:=-d\alpha\in\Omega_{T^{*}\ca}^{2}(T^*M)_{0}$. To see that $\omega$ is non-degenerated, consider its local expression:
$$\begin{array}{rcl}
 \qquad\qquad\omega= -d\alpha& = & -\left[\displaystyle\sum_{i}\sig[|q_i|]d(p_i)\wedge d( q_i )+\sum_i\sig[|q_i|+0\cdot 1+|p_i|\cdot 0]p_i\cdot d^2( q_i )\right]\qquad\qquad\\
	     & = & -\displaystyle\sum_{i}\sig[|q_i|+1\cdot 1+|q_i||p_i|]d( q_i )\wedge d(p_i)\\
	     & = & \displaystyle\sum_{i}d( q_i )\wedge d(p_i).  \hfill\square	   
\end{array}
$$
The non-degeneracy of $\omega$ allows us to give the following definition:
\begin{de}Let $(\ca,\omega)$ be a symplectic supermanifold. For any $f\in\om(M)$ let $X_f\in\tm(M)$ be the \textbf{Hamiltonian vector field} defined by the relation:
\[
 i(X_f)\omega=df,
\]
where $(i(X)\omega)(Y):=\sig[|X||\omega|]\omega(X,Y)=\omega(X,Y).$
\end{de}
Note that $X_f$ has necessarily the same parity as $f$ since $\omega$ is even. Direct calculation now shows:
\begin{prop}
 Let $\ca$ be a supermanifold and $(T^*\ca,\omega)$ be its cotangent bundle equipped with the canonical symplectic form. The Hamiltonian vector field $X_f\in\ca[T]_{T^*\ca}(T^*M)$ associated to an even function $f\in\ca[O]_{T^*\ca}(T^*M)_0$ is an even vector field and it is locally given by the following equations:
\[
 \begin{array}{ccl}
                      X_f(q_i)&=&\sig\pa[p_i]f\\
                      X_f(p_i)&=&-\pa[q_i]f.
                     \end{array}
\]
\end{prop}
\begin{cor}\label{XH}
In the situation of the preceding proposition, we have for $H$ the energy function the following equations: 
\[
 \begin{array}{ccl}
                    X_H(q_i)&=&\displaystyle  \sum_{j} p_j g^{ji}    \\
                     X_H(p_i)&=&  \displaystyle - \frac{1}{2}\sum_{k,j}\sig[|q_i||q_k|]p_k\pa[q_i](g^{kj}) p_j   .
                     \end{array}
\]
\end{cor}

\begin{de}Let $\ca$ be a supermanifold and $X$ be an even vector field on $\ca$. A \textbf{flow of} $X$ is a morphism $F:\bb\times\ca\to\ca$ of supermanifolds such that:
\begin{eqnarray}
  \widehat{\pa[t]}\circ F^*=F^*\circ X\text{ and }\:F\circ\text{inj}_{t=0}=\text{id}_{\ca},\label{eqflow}
\end{eqnarray}

where $t$ is the canonical coordinate on $\bb$, $\widehat{\pa[t]}$ the lift of the derivation $\pa[t]$ on $\bb$ to a derivation on $\bb\times\ca$ acting just on the first factor, and $\text{inj}_{t=0}:\ca\cong\{0\}\times\ca\overset{0\times \text{id}_{\ca}}{\longrightarrow}\bb\times\ca$ is the injection map at $t=0$.
\end{de}
\textbf{Remark.} Of course, the domain of the flow map $F$ is in general only an open subsupermanifold $\ca[V]$ of $\bb\times\ca$, containing $\{0\}\times\ca$. Abusively, we will often write $\bb\times\ca$ instead of $\ca[V]$.\\

The following result is due to Monterde and Montesinos \cite{Monterde}:
\begin{theo}
Let $X$ be an even vector field on a supermanifold $\ca$, $\widetilde{X}$ the reduced vector field on $M$, and $V_{\widetilde{X}}\subset\bb\times M$ the flow domain of $\widetilde{X}$. Then there exists a unique map $F:\ca[V]_X\to\ca$, where $\ca[V]_X$ is the subsupermanifold of $\bb\times\ca$ having as body $V_{\widetilde{X}}$, satisfying (\ref{eqflow}).

 \end{theo}
Given the flow $F$ of an even vector field $X$, we want to use an appropriate notion of initial condition to get ``integral supercurves". If we choose a point $q_0$ in $M$, we get a curve as follows:
\[
\begin{diagram}
 \node{\bb\cong\bb\times\{q_0\}}\arrow{e,t}{\text{inj}_{q=q_0}}\arrow{se,b}{\Psi}\node{\bb\times\ca}\arrow{s,r}{F}\\
 \node[2]{\ca}
\end{diagram}
\]
This type of initial condition leads to a poor notion of ``integral supercurves", since $\Psi$ factorizes over $M$, the body of the supermanifold $\ca=(M,\om)$. Thus we need a notion of initial condition that is more elaborate than simply choosing a point in the body of $\ca$.\\

Following Dumitrescu (\cite{dumi}) we give

\begin{de} Let $\ca$, $\ca[S]$ be supermanifolds, $X$ an even vector field on $\ca$ and $\ca[I]\in\text{Mor}(\ca[S],\ca)$. A morphism $\Phi:\bb\times\ca[S]\to\ca$ such that 
\begin{eqnarray}
  \widehat{\pa[t]}\circ \Psi^*=\Psi^*\circ X\text{ and }\:\Psi\circ\text{inj}_{t=0}=\ca[I],\label{eqintegracurv}
\end{eqnarray}
where $t$ is the canonical even coordinate on $\bb$ and $\text{inj}_{t=0}:\ca[S]\longrightarrow\bb\times\ca[S]$ the evaluation map at $t=0$, is called an ``\textbf{integral curve} for $X$, parametrized by $\ca[S]$ and satisfying the initial condition $\ca[I]$".
\end{de}
The case we are mostly concerned with will be $\ca[S]=\sbb[0|1]$ but we need this more general ``family notion" for integral curves notably in the proof of Theorem \ref{theoiso} below. We note that one can give the following geometrical description of the set of initial conditions (for $\ca[S]=\sbb[0|1]$):
\begin{prop}
 Let $\ca=(M,\om)$ be a supermanifold and $E$ be a real vector bundle over $M$ such that $\ca\cong(M,\Gamma^{\infty}_{\wedge E^*})$. Then there is a 1:1 correspondence between $E$ and $\text{Mor}(\sbb[0|1],\ca)$.
 \end{prop}
\textbf{Idea of the proof.} The local coordinates given by local trivializations of the vector bundle $E$ yield immediately
a bijection between $E$ and $\text{Mor}(\sbb[0|1],\ca)$.  \hfill $\square$\\

Similarly to the standard ungraded case, integral curves are given by factorization over the flow map:
\begin{prop}\label{uniqintcurv}(\cite{dumi})
 Let $X$ be an even vector field on $\ca$ and $F$ be the flow of $X$. Let $\ca[I]\in\text{Mor}(\ca[S],\ca)$ be an initial condition. There exists an unique integral curve $\Phi:\bb\times\ca[S]\to\ca$ of $X$ with initial condition $\ca[I]$. Precisely, it is given by $\Phi=F\circ(\text{id}_{\bb}\times\ca[I])$.
 \end{prop}
 We can now give the first fundamental results of this article, relating the integral curves of the geodesic flow on $T^*\ca$ to the geodesics of a Riemannian supermanifold $\ca$.
\begin{theo}\label{thetheorem}
 Let $(\ca,g)$ be a Riemannian supermanifold and $(T^{*}\ca,\omega,H)$ be the Hamiltonian dynamical system associated to $\ca$, and let $\ca[S]$ be a supermanifold. Furthermore, let $X_H$ be the Hamiltonian vector field associated to $H$, $\ca[I]$ an initial condition in $\text{Mor}(\ca[S],T^*\ca)$ and $\Psi$ the associated integral curve of $X_H$ ($\Psi=F\times(\text{id}_{\bb}\times\ca[I])$). Then the supercurve $\Phi:=\pi^{T^{*}\ca}\circ\Psi:\bb\times\ca[S]\to\ca$ satisfies the geodesic equation:
\[
 \frac{\nabla}{dt}\big(T\Phi^*(\pa[t])\big)=0,
\]
subject to the initial condition $\pa[t]\Phi\circ\text{inj}_{t=0}=g^{\flat}\circ\ca[I]\in\text{Mor}(\ca[S],T\ca)$, where $g^{\flat}:T^*\ca\overset{\cong}{\longrightarrow}T\ca$ is the isomorphism induced by the metric $g$.
\end{theo}
\textbf{Proof.} Let $\Psi$ be an integral curve of $X_H$. Then $\Psi$ is given by $\Psi= F\circ(\text{id}_{\bb}\times\ca[I])$ using Proposition (\ref{uniqintcurv}) where $ F:\bb\times T^{*}\ca\to T^{*}\ca$ is the flow of $X_H$ and $\ca[I]\in\text{Mor}(\ca[S],T^{*}\ca)$ is the initial condition.\\

Recall that $ F$ is the unique morphism such that 
\[ 
\widehat{\pa[t]}\circ  F^*= F^*\circ X_H\text{ and } F\circ\text{inj}_{t=0}=\text{id}_{T^{*}\ca}.
\]
In local coordinates $ F$ is uniquely determined by the following system of equations (compare Corollary \ref{XH}): 
\begin{eqnarray*}
 \widehat{\pa[t]}\left( F^{*}(q_i)\right)&=& \displaystyle\sum_{j}  F^*(p_jg^{ji})\\
 \widehat{\pa[t]}\left( F^{*}(p_i)\right)&=& \displaystyle- \frac{1}{2}\sum_{k,j}\sig[|i||k|] F^{*}\left(p_k\pa[q_i](g^{kj}) p_j\right)\\
(\text{inj}_{t=0})^{*}\circ F^*(q_i)&=& q_i\\
(\text{inj}_{t=0})^{*}\circ F^*(p_i)&=& p_i
\end{eqnarray*}
We note here and in the sequel the parity $|p_i|=|q_i|$ often by $|i|$.\\

Applying $(\text{id}_{\bb}\times\ca[I])^*$ to the two first equations of the above system and $\ca[I]^*$ to the remaning equations, we get
\begin{eqnarray}
\pa[t]\left(\Psi^{*}(q_i)\right)&=& \displaystyle\sum_{j} \Psi^*(p_jg^{ji})\label{sys1}\\
\pa[t]\left(\Psi^{*}(p_i)\right)&=& \displaystyle- \frac{1}{2}\sum_{k,j}\sig[|i||k|]\Psi^{*}\left(p_k\pa[q_i](g^{kj}) p_j\right)\label{sys2}\\
(\text{inj}_{t=0})^{*}\circ\Psi^*(q_i)&=& \ca[I]^*(q_i)\label{sys3}\\
(\text{inj}_{t=0})^{*}\circ\Psi^*(p_i)&=& \ca[I]^*(p_i)\label{sys4}
\end{eqnarray}
Equation (\ref{sys1}) gives now
\begin{equation}
\displaystyle\sum_{i} \pa[t]\Psi^{*}(q_i)\cdot\Psi^*(g_{ij})= \Psi^*(p_j). \label{sys1_2}
\end{equation}
Derivation along the ``even time direction" $t$ yields:
\begin{eqnarray*}
 \pa[t] \Psi^*(p_j) & = & \displaystyle\sum_{i}\left[ \pa[t]^2\Psi^{*}(q_i)\cdot\Psi^*(g_{ij})+\pa[t]\Psi^{*}(q_i)\cdot\pa[t]\Psi^*(g_{ij})\right]\\
		    & = & \displaystyle\sum_{i} \pa[t]^2\Psi^{*}(q_i)\cdot\Psi^*(g_{ij})+\sum_{k,i}\pa[t]\Psi^{*}(q_i)\cdot\pa[t]\Psi^*(q_k)\cdot\Psi^*(\paq[k] g_{ij}).
\end{eqnarray*}
Let us first consider the last summand of the last RHS:
\[
 \begin{array}{l}
\displaystyle\sum_{k,i}\pa[t]\Psi^{*}(q_i)\cdot\pa[t]\Psi^*(q_k)\cdot\Psi^*(\paq[k] g_{ij})\\
\hspace{1cm}\vspace*{2mm}=\displaystyle\frac{1}{2}\left[\sum_{k,i}\pa[t]\Psi^{*}(q_i)\cdot\pa[t]\Psi^*(q_k)\cdot\Psi^*(\paq[k] g_{ij}) +\sum_{i,k}\pa[t]\Psi^{*}(q_k)\cdot\pa[t]\Psi^*(q_i)\cdot\Psi^*(\paq[i] g_{kj})\right]\\
\hspace{2cm}=\displaystyle\frac{1}{2}\sum_{k,i}\pa[t]\Psi^{*}(q_i)\cdot\pa[t]\Psi^*(q_k)\left(\Psi^*(\paq[k] g_{ij})+\sig[|i||k|]\Psi^*(\paq[i] g_{kj})  \right) .
\end{array}
\]

After relabelling indices we arrive at
\begin{equation}
 \pa[t] \Psi^*(p_l) = \displaystyle\sum_{j} \pa[t]^2\Psi^{*}(q_j)\cdot\Psi^*(g_{jl})+\frac{1}{2}\sum_{i,j}\pa[t]\Psi^{*}(q_j)\cdot\pa[t]\Psi^*(q_i)\cdot\left(\Psi^*(\paq[i] g_{jl})+\sig[|j||i|]\Psi^*(\paq[j] g_{il}) \right).\label{sys1_3}
\end{equation}
Using equation (\ref{sys1_2}), equation (\ref{sys2}) yields:
\begin{eqnarray*}
 \pa[t]\left(\Psi^{*}(p_i)\right)&=& \displaystyle- \frac{1}{2}\sum_{k,j}\sig[|i||k|]   \Psi^{*}(p_k)  \cdot  \Psi^{*}(\pa[q_i]g^{kj})\cdot   \Psi^{*}( p_j)   \\
				 &=& \displaystyle- \frac{1}{2}\sum_{\substack{ k,j \\ l,s }}\sig[|i||k|]     \pa[t]\Psi^{*}(q_l)\cdot\Psi^*(g_{lk})  \cdot    \Psi^{*}(\pa[q_i]g^{kj})\cdot    \pa[t]\Psi^{*}(q_s)\cdot\Psi^*(g_{sj}).    \\ 
\end{eqnarray*}
Summing the following equality over $k$
\[
\pa[q_i]( g_{lk}\cdot  g^{kj} ) = \pa[q_i]g_{lk}\cdot g^{kj} + \sig[|i|(|l|+|k|)]g_{lk}\cdot\pa[q_i]g^{kj}
\]
gives
\[
 \displaystyle\sum_k \sig[|i||l|] \pa[q_i]g_{lk}\cdot g^{kj} + \sum_k\sig[|i||k|]g_{lk}\cdot\pa[q_i]g^{kj} = 0.
\]
Then
\begin{eqnarray*}
 \pa[t]\left(\Psi^{*}(p_i)\right)&=&  \displaystyle \frac{1}{2}\sum_{\substack{ k,j \\ l,s }}\sig[|i||l|]     \pa[t]\Psi^{*}(q_l)\cdot\Psi^*(\pa[q_i]g_{lk})\cdot      \Psi^{*}(g^{kj})\cdot   \pa[t]\Psi^{*}(q_s)\cdot\Psi^*(g_{sj})    \\ 
				  &=&  \displaystyle \frac{1}{2}\sum_{\substack{ k,j \\ l,s }}\sig[|i||l|+|s|]     \pa[t]\Psi^{*}(q_l)\cdot\Psi^*(\pa[q_i]g_{lk})\cdot      \Psi^{*}(g^{kj})   \Psi^*(g_{js})\cdot\pa[t]\Psi^{*}(q_s)    \\ 
				  &=&  \displaystyle \frac{1}{2}\sum_{ k,l}\sig[|i||l|+|k|]     \pa[t]\Psi^{*}(q_l)\cdot\Psi^*(\pa[q_i]g_{lk})   \cdot   \pa[t]\Psi^{*}(q_k)    \\ 
				   &=&  \displaystyle \frac{1}{2}\sum_{ k,l}\sig[|k|(|i|+|l|)+|i||l|]     \pa[t]\Psi^{*}(q_l)\cdot \pa[t]\Psi^{*}(q_k)\cdot\Psi^*(\pa[q_i]g_{lk})      .\\ 
\end{eqnarray*}
After again relabelling indices, we have
\begin{eqnarray*}
 \pa[t]\Psi^{*}(p_l)= \displaystyle\frac{1}{2}\sum_{ j,i}\sig[|j|(|l|+|i|)+|l||i|]     \pa[t]\Psi^{*}(q_i)\cdot  \pa[t]\Psi^{*}(q_j) \cdot \Psi^*(\pa[q_l]g_{ij})
 \end{eqnarray*}
 \begin{equation}
  \hspace{.3cm}
  = \displaystyle\frac{1}{2}\sum_{ j,i}\sig[|l|(|i|+|j|)]  \pa[t]\Psi^{*}(q_j) \cdot   \pa[t]\Psi^{*}(q_i) \cdot  \Psi^*(\pa[q_l]g_{ij})  .    \label{sys2_2}
 \end{equation}
Multiplying (\ref{sys1_3}) with $\Psi^*(g^{lk})$ and summing over $l$ gives: 
\begin{equation}
  \displaystyle \sum_l \pa[t] \Psi^*(p_l)\cdot\Psi^*(g^{lk}) = \pa[t]^2\Psi^{*}(q_k)+\frac{1}{2}\sum_{i,j , l}\pa[t]\Psi^*(q_j)\cdot\pa[t]\Psi^{*}(q_i)\cdot\Psi^*\left(\paq[i] g_{jl}+\sig[|j||i|]\paq[j] g_{il} \right)\cdot\Psi^*(g^{lk}).\label{sys1_4}
\end{equation}
Applying the same operations to (\ref{sys2_2}) yields
\begin{equation}
 \displaystyle \sum_l \pa[t] \Psi^*(p_l)\cdot\Psi^*(g^{lk}) =  \displaystyle \frac{1}{2}\sum_{j,i,l }\sig[|l|(|i|+|j|)]     \pa[t]\Psi^{*}(q_j) \cdot \pa[t]\Psi^{*}(q_i)\cdot \Psi^*(\pa[q_l]g_{ij})\cdot \Psi^*(g^{lk})   .    \label{sys2_3}
\end{equation}

Now, using the equality $\Phi^*(q_i)=\Psi^*\circ(\pi^{T^*\ca})^*(q_i)=\Psi^*(q_i)$ and equating the RHS of the equalities (\ref{sys1_4}) and (\ref{sys2_3}), we get exactly (\ref{equageodesic}) since the Christoffel symbols are given by (\ref{chritoffel}).\\

Let us check that $\pa[t]\Phi\circ\text{inj}_{t=0}=g^{\flat}\circ\ca[I]$.\\

Recall the canonical sections $s_0,s_1$ defined before Theorem \ref{uniqgeo}. Since $\Phi=\pi^{T^*\ca}\circ F\circ(\text{id}_{\bb}\times\ca[I])$, we can, in fact, show that $\pa[t]\Phi=g^{\flat}\circ F\circ(\text{id}_{\bb}\times\ca[I])$.\\
Firstly, we have
$$
\begin{array}{rcl}
\pa[t]\Phi & = & T\Phi\circ(s_1\times s_0) \\
		   & = & T\pi^{T^*\ca}\circ T F\circ(\text{id}_{T\bb}\times T\ca[I])\circ(s_1\times s_0)\\
		   & = & T\pi^{T^*\ca}\circ T F\circ(s_1\times( T\ca[I]\circ s_0)):\bb\times\ca[S]\to T\ca.\\
\end{array}
$$
If $(q_i,p_i,v_i,r_i)$ are the usual local coordinates on $T(T^*\ca)$ and $(q_i,v_i)$ local coordinates on $T\ca$, identify with the image of the tangent map associated to the zero section of $T^*\ca\to\ca$. Let further denote $(v_t,t)$ the point $v_t\frac{d}{dt}_{|t}$ of $T\bb$. We then have
$$
\begin{array}{rcl}
(\pa[t]\Phi)^*(q_i) & = & (s_1^*\times( s_0^*\circ(T\ca[I])^*)) \circ (T F)^* \circ (T\pi^{T^*\ca})^*(q_i)\\
                    & = & (s_1^*\times( s_0^*\circ(T\ca[I])^*)) \circ (T F)^*  (q_i)\\
                    & = & (s_1^*\times( s_0^*\circ(T\ca[I])^*))(  F^*  (q_i))\\
                    & = & (\text{id}_{\bb}^*\times \ca[I]^*)(  F^*  (q_i)) \\
                    & = & (\text{id}_{\bb}^*\times \ca[I]^*)\circ  F^* \circ (g^\flat)^*(q_i) \\
                    & = & \left(g^\flat\circ F\circ(\text{id}_{\bb}\times \ca[I])\right)^*(q_i),
\end{array}
$$
and
$$
\begin{array}{rcl}
(\pa[t]\Phi)^*(v_i) & = & (s_1^*\times( s_0^*\circ(T\ca[I])^*)) \circ (T F)^* \circ (T\pi^{T^*\ca})^*(v_i)\\
                    & = & (s_1^*\times( s_0^*\circ(T\ca[I])^*)) \circ (T F)^*  (v_i)\\
                    & = & (s_1^*\times( s_0^*\circ(T\ca[I])^*))\left(\sum_{j}\left(v_j.\pa[q_j] F^*(q_i)+r_j.\pa[p_j] F^*(q_i)\right)+v_t.\pa[t] F^*(q_i)\right)\\
                    & = & (\text{id}_{\bb}^*\times \ca[I]^*)\left( \pa[t] F^*(q_i) \right) \\
                    & = & (\text{id}_{\bb}^*\times \ca[I]^*)\left(( F^*\circ X_H)(q_i)\right),\quad\text{since } F\text{ is the flow of }X_H, \\
                    & = & (\text{id}_{\bb}^*\times \ca[I]^*)\left( F^*(\sum_{j}p_j g^{ji})\right) \\
                    & = & (\text{id}_{\bb}^*\times \ca[I]^*)\left( F^*\circ(g^\flat)^*(v_i)\right) \\
                    & = & \left(g^\flat\circ F\circ(\text{id}_{\bb}\times \ca[I])\right)^*(v_i).
\end{array}
$$
Since we now have $\pa[t]\Phi=g^\flat\circ F\circ(\text{id}_{\bb}\times\ca[I])$, we can easily conclude the proof of the theorem:
\begin{equation}
\pa[t]\Phi\circ\text{inj}_{t=0}=g^{\flat}\circ F\circ(\text{id}_{\bb}\times\ca[I])\circ\text{inj}_{t=0}=g^{\flat}\circ\ca[I].\tag*{$\square$}
\end{equation}
Reciprocally, we have the following:
\begin{theo}\label{rthetheorem}
 Let $\Phi:\bb\times\ca[S]\to\ca$ be a supergeodesic with initial condition $\pa[t]\Phi\circ\text{inj}_{t=0}=\ca[J]$ in $\text{Mor}(\ca[S],T\ca)$. Then  $g^{\sharp}\circ \pa[t]\Phi:\bb\times\ca[S]\to T^{*}\ca$ is an integral curve of $X_H$ with initial condition $g^{\sharp}\circ\ca[J]\in\text{Mor}(\ca[S],T^*\ca)$, where $g^{\sharp}:T\ca\overset{\cong}{\longrightarrow} T^*\ca$ is the isomorphism induced by the metric $g$.
\end{theo}
\textbf{Proof.} Let $\Phi:\bb\times\ca[S]\to\ca$ be a supercurve such that $\frac{\nabla}{dt}\big(T\Phi^*(\pa[t])\big)=0$ and $\pa[t]\Phi\circ\text{inj}_{t=0}=\ca[J]$ in $\text{Mor}(\ca[S],T\ca)$. We define $\ca[I]:=g^\sharp\circ\ca[J]$ and $\Psi:= F\circ(\text{id}_{\bb}\times \ca[I])$.
By the previous theorem, we know that $\pi^{T^*\ca}\circ\Psi$ is a geodesic of $\ca$ with initial condition $g^\flat\circ\ca[I]=\ca[J]$. By unicity, we have $\Phi=\pi^{T^*\ca}\circ\Psi$. Now, using a calculation made in the proof of the preceding theorem, we have $\pa[t]\Phi=g^\flat\circ F\circ(\text{id}_{\bb}\times\ca[I])$.
Then $g^\sharp\circ\pa[t]\Phi= F\circ(\text{id}_{\bb}\times\ca[I])$, i.e., $g^\sharp\circ\pa[t]\Phi$ is the integral curve of $X_H$ with initial condition $\ca[I]=g^\sharp\circ\ca[J]$.\hfill $\square$

\section{The exponential map of a Riemannian supermanifold}
In this section we study, for $q$ in the body of a supermanifold $\ca$, the fiber $T_{q}\ca$ of the tangent bundle over
$q$, the restriction of the tangent map of a morphism to this fiber, and its relation to the ``numerical tangent map"
(Lemma \ref{maptnum}). The main result of this section is that a Riemannian supermanifold possesses a (Riemannian) exponential
map and that $\text{exp}_q:T_q\ca\to\ca$ is a diffeomorphism near $0$ (Thm. \ref{expdiffeo}). As an example of an application we show that linearization in a fixed point of an isometry of a Riemannian supermanifold is faithful (Cor. \ref{faithful}).\\

Let $(\ca,g)$ be a Riemannian supermanifold, and $X_H$ be the even Hamiltonian vector field on $T^*\ca$ associated to the energy function $H$.
Let furthermore $F:\bb\times T^*\ca[M]\to T^*\ca[M]$ be the flow of $X_H$ and $g^\sharp:T\ca\to T^*\ca$ be the isomorphism of supermanifolds induced by the metric $g$.\\

Consider the following composition of applications:\\
$$\begin{diagram}
\node{\bb\times T\ca}\arrow{e,t}{id_{\bb} \times g^\sharp}\node{\bb\times T^*\ca}\arrow{e,t}{F}\node{T^*\ca}\arrow{e,t}{\pi^{T^*\ca}}\node{\ca.}
\end{diagram}$$
We want to take the fibre of $T\ca$ above $q\in M$:
$$\begin{diagram}
\node[2]{T\ca}\arrow{s,r}{\pi^{T\ca}}\\
\node{\{q\}}\arrow{e,t}{q}\node{\ca}
\end{diagram}$$
Let us recall the morphism of rings of global sections $(\pi^{T\ca})^*(TM):\om(M)\to\ca[O]_{T\ca}(TM)$, associated to the sheaf morphism $(\pi^{T\ca})^*:(\widetilde{\pi^{T\ca}})^{-1}(\om)\to\ca[O]_{T\ca}$ over $TM$.\\

Let $\ca[J]_q$ be the ideal of $\om(M)$ defined by $$\ca[J]_{q}:=\{f\in\om(M) \: | \: q^*(f)=0\},$$ where $q^*(f)=\widetilde{f}(q)$.
We denote the ideal in $\ca[O]_{T\ca}(TM)$ generated by $(\pi^{T\ca})^*(\ca[J]_{q})$ by 
$$\ca[J]^{\pi^{T\ca}}_{q}:=\langle(\pi^{T\ca})^*(\ca[J]_{q}))\rangle_{\ca[O]_{T\ca}(TM)}.$$
Then (compare \cite{Leites} for this construction in general), there exists a subsupermanifold $\ca[S]\overset{i_q}{\hookrightarrow}T\ca$ of $T\ca$ such that
$$
\ca[J]_{\ca[S]}:=\{f\in\ca[O]_{T\ca}(TM) \:|\: i_q^*(f)=0\}=\ca[J]^{\pi^{T\ca}}_{q}.
$$
The supermanifold $\ca[S]$ will be denoted by $T_{q}\ca$.\\

Note that if $\ca=\sbb[m|n]$ and $q$ is an element of $\bb^m$, we have $T\ca=\sbb[m|n]\times\sbb[m|n]$ and $T_q\ca=\sbb[m|n]\times\{q\}\cong\sbb[m|n]$. In general, given a finite dimensional $\bb[Z]_2$-graded $\bb$-vector space $V=V_0\oplus V_1$, we have of course an associated ``linear supermanifold". Then $T_q\ca$ is the supermanifold associated to $T^\text{num}_q\ca$.\\

We now define the map $\text{exp}_q:T_q\ca\to\ca$ via the following diagram:
$$\begin{diagram}
\node{\bb\times T\ca}\arrow{e,t}{id_{\bb} \times g^\sharp}\node{\bb\times T^*\ca}\arrow{e,t}{F}\node{T^*\ca}\arrow{e,t}{\pi^{T^*\ca}}\node{\ca}\\
\node{\bb\times T_{q}\ca}\arrow{n,r}{id_{\bb} \times i_q}\\
\node{\{1\}\times T_{q}\ca\cong T_{q}\ca}\arrow{n,r}{1\times id_{T_{q}\ca}}\arrow{nneee,b}{\text{exp}_{q}}
\end{diagram}$$
Recall that the flow $F$ is only defined on an open subsupermanifold of $\bb\times T^*\ca$ containing $\{0\}\times T^*\ca$. Then $\text{exp}_q$ is in general only defined on an open subsupermanifold of $T_q\ca$ containing in its body $(0,q)$. We can now show that $\text{exp}_q$ is a locally diffeomorphism near $0\in T_qM$.
\begin{lem}
Let $\Phi:\ca\to \ca[N]$ be a morphism of supermanifolds. Let $q\in M$, then we have a unique morphism of supermanifolds $T_q\Phi:T_q\ca\to T_{\widetilde{\Phi}(q)}\ca[N]$ such that $i_{\widetilde{\Phi}(q)}\circ T_q\Phi=T\Phi\circ i_q$.

\end{lem}
\textbf{Proof.} We have the following commutative diagram:
\[
\begin{diagram}
\node{T_q\ca}\arrow{e,t}{i_q}\arrow{s}\node{T\ca}\arrow{e,t}{T\Phi}\arrow{s,r}{\pi^{T\ca}}\node{T\ca[N]}\arrow{s,r}{\pi^{T\ca[N]}}\node{T_{p}\ca[N]}\arrow{w,t}{i_p}\arrow{s}\\
\node{\{q\}}\arrow{e}\node{\ca}\arrow{e,t}{\Phi}\node{\ca[N]}\node{\{p\},}\arrow{w}
\end{diagram}
\]
where $p:=\widetilde{\Phi}(q)\in N$.
We will show that the morphism $T\Phi\circ i_q$ factorizes over the morphism $i_p$.\\

Recall that a morphism of real smooth supermanifolds is completely determined by the superalgebra morphism given by the pulling-back global sections.

Let us define the morphism $T_q\Phi:T_q\ca\to T_{p}\ca[N]$ by setting its body map equal to $\widetilde{T_q\Phi}:=T_q\widetilde{\Phi}$ and on the level of global sections by
$$\begin{array}{rcl}
(T_q\Phi)^*:\ca[O]_{T_p\ca[N]}(T_p N)&\longrightarrow&\ca[O]_{T_q\ca}(T_q M)\\
               					f     	   &\longmapsto    &i_q^*\circ (T\Phi)^*(\hat{f}\,),
\end{array}$$
where $\hat{f}$ is an element of $\ca[O]_{T\ca[N]}(TN)$ such that $i_p^*(\hat{f}\,)=f$.\\

For all $\hat{f}\in\ca[O]_{T\ca[N]}(TN)$ such that $i_p^*(\hat{f}\,)=0$, we have $\hat{f}\in\ca[J]_{T_p\can}=\ca[J]^{\pi^{T\ca[N]}}_{p}$ by the very construction of $T_{p}\ca[N]$ and then by definition of $\ca[J]^{\pi^{T\ca[N]}}_{p}$: 
$$\hat{f}\in\langle(\pi^{T\ca[N]})^*(\ca[J]_{p})\rangle_{\ca[O]_{T\ca[N]}(TN)}.$$
Thus it is enough to show that $i_q\circ(T\Phi)^*(\hat{f}\,)=0$ whenever $\hat{f}=(\pi^{T\can})^*(g)$ with $g\in\ca[J]_p$ in order to assure that $(T_q\Phi)^*$ is well-defined.\\
For $g\in\ca[J]_p\lhd\on(N)$ we have
$$
(T\Phi)^*\left((\pi^{T\can})^*(g)\right)=(\pi^{T\ca})^*(\Phi^*(g))\in\langle(\pi^{T\ca})^*(\ca[J]_q)\rangle_{\ca[O]_{T\ca}(TM)},
$$
since the identity $\widetilde{\Phi^*(h)}(q)=\widetilde{h}(p)$ for $h\in\on(N)$ shows that $g\in\ca[J]_p$ implies that $\Phi^*(g)\in\ca[J]_q$.\\
Moreover, since $\ca[J]_{T_q\ca}=\langle(\pi^{T\ca})^*(\ca[J]_q)\rangle_{\ca[O]_{T\ca}(TM)}$ we have $i_q^*\left((\pi^{T\ca})^*(\Phi^*(g))\right)=0$, i.e. $i_q^*\circ(T\Phi)^*(\hat{f}\,)=0$.   
\hfill$\square$\\

Let us recall that for a morphism of supermanifolds $\Phi:\ca\to \ca[N]$, and a point $q\in M$, we have a numerical analogue of the tangent map, i.e. a linear map $T^{\text{num}}_q\Phi:T_q^{\text{num}}\ca\to T_{\widetilde{\Phi}(q)}^{\text{num}}\ca[N]$ defined by 
\[
  T^{\text{num}}_q\Phi(X)\big([f]_{\widetilde{\Phi}(q)} \big):= X\big([\Phi^*(f)]_q\big),\quad \forall X\in T^{\text{num}}_q\ca,\quad\forall  [f]_{\widetilde{\Phi}(q)}\in \ca[O]_{\ca[N],\widetilde{\Phi}(q)}.
\]
 It turns out that this map encodes exactly the same information as  the morphism $T_q\Phi$.
\begin{lem}\label{maptnum}
Let $\Phi:\ca\to\can$ be a morphism of supermanifolds and $q\in M$. Then the following holds:
\begin{enumerate}[(i)]
\item $T_q\Phi:T_q\ca\to T_{\widetilde{\Phi}(q)}\can$ is the morphism of (linear) supermanifolds associated to the even $\bb$-linear map $T^{\text{num}}_q\Phi:T^{\text{num}}_q\ca\to T^{\text{num}}_{\widetilde{\Phi}(q)}\can$, and
\item the numerical tangent map $T^{\text{num}}_q\Phi$ is uniquely determined by the morphism of supermanifolds $T_q\Phi$.
\end{enumerate}
\end{lem}
\textbf{Proof. }A calculation in the local coordinates on $T\ca$ and $T\can$, that are induced from local coordinates on $\ca$ and $\can$, directly shows that both maps $T_q\Phi$ and $T^{\text{num}}_q\Phi$ are uniquely determined by the evaluation of the super Jacobi matrix of $\Phi$ in the point $q$ of the body of $\ca$.\hfill$\square$\\

It is well-known that a morphism $\Phi:\ca\to\ca[N]$ is a local diffeomorphism near $q$ if and only if $T^{\text{num}}_q\Phi$ is an isomorphism of $\bb[Z]_2$-graded vector space. (see, e.g., \cite{var}, Theorem 4.4.1)\\

Using the preceding lemma, we thus immediately get:
\begin{prop}
The morphism $\Phi:\ca\to\ca[N]$ is locally a diffeomorphism near $q\in M$ if and only if $\:T_q\Phi:T_q\ca\to T_{\widetilde{\Phi}(q)}\ca[N]$ is an isomorphism of supermanifolds. 
\end{prop}
\begin{theo}\label{expdiffeo}
The map $\text{exp}_q:T_q\ca\to\ca$ is a local diffeomorphism near $0$.
\end{theo}
\textbf{Proof.} We will show that the map $T_0(\text{exp}_q):T_0(T_q\ca)\to T_q\ca$ is an isomorphism. More precisely, after identifying $T_0(T_q\ca)$ with $T_q\ca$ as in the classical case, we will see that $T_0(\text{exp}_q)=\text{id}_{T_q\ca}$.\\

Let us focus on the map:
$$\begin{array}{rcl}
T^{\text{num}}_0(\text{exp}_q):T^{\text{num}}_0(T_q\ca)=\text{Der}\left((\ca[O]_{T_q\ca,0}),\bb\right)&\longrightarrow&\text{Der}\left((\ca[O]_{\ca,q}),\bb\right)=T^{\text{num}}_q\ca\\
               					X     	   &\longmapsto    &X\circ \text{exp}_q^*
\end{array}.$$
If $(q_i)$ is a system of local coordinates on $\ca$ around $q$, we denote $(q_i,v_i)$ the induced local coordinates on $T\ca$. Let us denote the natural coordinates on $T_q\ca$ by $(v_i)$ as well.
Then $T^{\text{num}}_0(T_q\ca)$ is the $\bb[Z]_2$-graded $\bb$-vector space spanned by $(\partial_{v_i})$ where $|\partial_{v_i}|=|v_i|$.\\

The classical theory of exponential maps on ungraded Riemannian manifolds gives the following equality
$$
T^{\text{num}}_0(\text{exp}_q)(\partial_{v_i})=\partial_{q_i}\quad\forall i\text{ such that }|\partial_{v_i}|=0.
$$
Let us now deal with the odd case. If $X=\partial_{v_i}$ with $|v_i|=1$, the ``integral curve" $\gamma:\sbb[0|1]\to T_q\ca$ of $X$ is defined by $\gamma^*(v_j)=\delta_{i,j}\cdot \theta$, where $\theta$ is the canonical coordinate of $\sbb[0|1]$.\\

Since
$$
X(f)=\frac{d}{d\theta}\gamma^*(f)\quad\forall f\in\ca[O]_{T_q\ca},
$$ 
we have 
$$
T^{\text{num}}_0(\text{exp}_q)(X)=\frac{d}{d\theta}\gamma^*\circ \text{exp}_q^*,
$$
where $\gamma^*\circ \text{exp}_q^*=\gamma^*\circ\left(\text{inj}^*_{t=1}\times\left(i_q^*\circ g^{\sharp*} \right)\right)\circ F^*\circ \left(\pi^{T^*\ca}\right)^*=\left(\text{inj}^*_{t=1}\times\left(\gamma^*\circ i_q^*\circ g^{\sharp*} \right)\right)\circ F^*\circ \left(\pi^{T^*\ca}\right)^*$.\\
Let us consider the map $\Gamma:\sbb\to\ca$ defined by:
$$
\Gamma^*:=\left(\text{id}^*_{\bb}\times\left(\gamma^*\circ i_q^*\circ g^{\sharp*} \right)\right)\circ F^*\circ \left(\pi^{T^*\ca}\right)^*.
$$
Then, using Theorem \ref{thetheorem}, $\Gamma$ is the unique geodesic of $\ca$ with initial condition $\eta= i_q\circ\gamma\in\text{Mor}\left(\sbb[0|1],T\ca\right)$, i.e. $\Gamma$ is the supercurve such that
$$
\left\lbrace\begin{array}{rclr}
\displaystyle \pa[t]^2 \Gamma^*(q_k)  + \sum_{|u|=|l|=0}  \pa[t]\Gamma^*(q_u)\pa[t]\Gamma^*(q_l)\Gamma^*(\Gamma_{lu}^k) & = & 0 & \text{ if } |q_k|=0,\vspace{3mm}
\\
\displaystyle \pa[t]^2 \Gamma^*(q_k)  + \sum_{\substack{|u|=0 \\ |l|=0}}  \pa[t]\Gamma^*(q_u)\pa[t]\Gamma^*(q_l)\Gamma^*(\Gamma_{lu}^k)+2\sum_{\substack{|u|=0 \\ |l|=1}}  \pa[t]\Gamma^*(q_u)\pa[t]\Gamma^*(q_l)\Gamma^*(\Gamma_{lu}^k) & = & 0 & \text{ if } |q_k|=1,
\end{array}\right.
$$
satisfying the initial condition
$$
\begin{array}{cc}
\left\lbrace\begin{array}{rclc}
\Gamma^*(q_k)_{|t=0}&=&\widetilde{q_k}(q) & \text{ if } |q_k|=0,\\
\Gamma^*(q_k)_{|t=0}&=&0 &\text{ if } |q_k|=1 ,
\end{array}
\right. & 
\left\lbrace\begin{array}{rclc}
 \partial_t\Gamma^*(q_k)_{|t=0}&=&\eta^*(v_k)=0 &\text{ if } |q_k|=0 ,  \\
 \partial_t\Gamma^*(q_k)_{|t=0}&=&\eta^*(v_k)=\delta_{k,i}\cdot \theta & \text{ if } |q_k|=1.
\end{array}
\right.
\end{array}
$$
Then the unique solution $\Gamma:\sbb\to\ca$ of this system is 
\[
\left\lbrace\begin{array}{rclc}
\Gamma^*(q_k)&=&\widetilde{q_k}(q) & \text{ if }|q_k|=0\\
\Gamma^*(q_k)&=&t\cdot \theta\cdot \delta_{k,i} &\text{ if } |q_k|=1 .
\end{array}
\right. 
\]
Now, we have
$$
T^{\text{num}}_0(\text{exp}_q)(X)=\frac{d}{d\theta}\text{inj}_{t=1}^*\circ\Gamma^*\in\text{Der}\left(\ca[O]_{\ca,q},\bb\right),
$$
and then
$$
T^{\text{num}}_0(\text{exp}_q)(X)(q_k)=\frac{d}{d\theta}\Gamma^*(q_k)_{|t=1}=\left\lbrace
\begin{array}{cl}
0 &\text{ if }|q_k|=0 \\
\delta_{k,i}&\text{ if }|q_k|=1,
\end{array}\right.
$$
that is 
$$
T^{\text{num}}_0(\text{exp}_q)(X)=\partial_{q_i}.
$$
This finishes the proof that $T_0(\text{exp}_q)$ is the identity map. Therefore $\text{exp}_q$ is a local diffeomorphism near $0$.\hfill$\square$\\

\textbf{Remark.} Already calculations with explicit Riemannian metrics on $\ca = {\mathbb{R}^{1|2}}$ show that the
most direct extension of the ungraded Gauss lemma can not hold true. I.e., given a Riemannian supermanifold $\ca$, then, in general,
for $q\in M$, $v=v_0+v_1\in T^{\text{num}}_q\ca$, and $w\in T_{v_0}^{\text{num}}(T_{q}\ca)\cong T_q^{\text{num}}\ca$, 
$$
g_{\widetilde{\text{exp}}_q(v_0)}\left((T^{\text{num}}_{v_0}\text{exp}_q)(v),(T^{\text{num}}_{v_0}\text{exp}_q)(w)\right)\ne g_{q}\left(v,w\right).
$$
\\

We can now apply the properties of the above constructed super exponential map to generalize a fundamental linearization result on Riemannian isometries to the category of supermanifolds.\\

Let $(\ca,g^{\ca})$ and $(\ca[N],g^{\ca[N]})$ be Riemannian supermanifolds and $\Phi:\ca\to\ca[N]$ a morphism.\\
Recall that $g^{\ca[N]}$ can be viewed as an element of $\Omega^1_{\can}\otimes_{\on} \Omega^1_{\can}$, moreover $\Phi^*:\widetilde{\Phi}^{-1}\on\to\om$ induces a morphism of sheaves
$$
\begin{diagram}
\node{\widetilde{\Phi}^{-1}\Omega^1_{\can}}\arrow{se}\arrow[2]{e,t}{\Phi^*}\node[2]{\Omega^1_{\ca}}\arrow{sw}\\
\node[2]{M}
\end{diagram}
$$
Thus the expression $\Phi^*(\widetilde{\Phi}^{-1}g^{\can})$ makes sense as an element of $\Omega^1_{\ca}\otimes_{\om}\Omega^1_{\ca}$.
\begin{de}
Let $(\ca,g^{\ca})$ and $(\ca[N],g^{\ca[N]})$ be riemannian supermanifolds. An \textbf{isometry} is an
isomorphism of supermanifolds $\Phi:\ca\to\can$ such that the following equality is satisfied
$$
\Phi^*(\widetilde{\Phi}^{-1}g^{\can})=g^{\ca}.
$$
\end{de}
If $(q_i^{\ca})$ and $(q_j^{\can})$ are local coordinates on $\ca$ and $\can$ respectively, the above condition means that 
$$
g^{\ca}_{ij}=\displaystyle \sum_{k,l}\sig[|q_k^{\can}|(|q_j^{\ca}|+|q_l^{\can}|)]\partial_i^{\ca}(\Phi^*(q_k^{\can}))\cdot\partial_j^{\ca}(\Phi^*(q_l^{\can}))\cdot\Phi^*\left(g^{\can}_{kl}\right),
$$
where $|\partial_i^{\ca}|$ and $|\partial_j^{\can}|$ denote $|\partial_{q_i^{\ca}}|$ and $|\partial_{q_j^{\can}}|$ respectively.
\begin{theo}\label{theoiso}
Let $\Phi:\ca\to\can$ be an isometry and $q\in M$, then we have 
$$
\displaystyle\Phi\circ\text{exp}\,_q^{\ca}=\text{exp}\,_{\widetilde{\Phi}(q)}^{\can}\circ T_q\Phi.
$$
\end{theo}
\textbf{Proof.} Let us define the following supercurves 
$$
\Lambda:=\pi^{T^*\ca}\circ F^{\ca}\circ(\mathbb{1}_{\bb}\times(g^{\ca\sharp}\circ i_q)):\bb\times T_q\ca\to\ca
$$
and
\begin{equation}\label{curvegamma}
\Gamma:=\Phi\circ\Lambda:\bb\times T_q\ca\to\can,
\end{equation}
where $F^{\ca}$ is the geodesic flow of the Hamiltonian dynamical system $(T^*\ca,\omega,H)$.\\

We want to show that $\Gamma$ is a geodesic of $\can$ with initial condition $\eta=i_{\widetilde{\Phi}(q)}\circ T_q\Phi\in\text{Mor}(T_q\ca, T\can)$. By the unicity of geodesics, it will then follow 
\begin{equation}
\Gamma=\pi^{T^*\can}\circ F^{\can}\circ(\mathbb{1}_{\bb}\times(g^{\can\sharp}\circ i_{\widetilde{\Phi}(q)}\circ T_q\Phi)).\label{gamma}
\end{equation}
The statement of the theorem will then be obtained by evaluating both expressions (\ref{curvegamma}) and (\ref{gamma}) for $\Gamma$ at $t=1$.\\

Since $\Phi$ is an isometry and the Levi-Civita connection is unique
$$
\left(\Phi^*\nabla^{\can}\right)(T\Phi^*)=T\Phi^*(\nabla^{\ca})
$$
and thus
$$
\frac{\nabla^{\can}}{dt}\left(T\Gamma^*(\pa[t])\right)=T\Phi^*\left(\frac{\nabla^{\ca}}{dt}\left(T\Lambda^*(\pa[t])\right)\right)=0
$$
because $\Lambda$ is a geodesic of $\ca$. It follows that $\Gamma$ is a geodesic.\\

Furthermore
$$
\begin{array}{ccl}
\pa[t]\Gamma\circ\text{inj}_{t=0} & = & \pa[t](\Phi\circ\Lambda)\circ\text{inj}_{t=0} \\                    
				    & = & T\Phi\circ\pa[t]\Lambda\circ\text{inj}_{t=0}\\
                    & = & T\Phi\circ i_q\\
                    & = & i_{\widetilde{\Phi}(q)}\circ T_q\Phi=\eta.						     
\end{array}
$$
We thus have shown that (\ref{gamma}) holds and the claim of the theorem follows by evaluating both sides in $t=1$.\hfill $\square$\\

Now, we can give the announced linearization result.
\begin{cor}\label{faithful}
Let $\ca$ be a connected Riemannian supermanifold and let $\Phi:\ca\to\ca$ be an isometry such that 
$$
\widetilde{\Phi}(q_0)=q_0\text{ and }T_{q_0}\Phi=\text{id}_{T_{q_0}\ca}\text{ for a fixed point }q_0\in M.
$$
Then $\Phi$ is the identity morphism of $\ca$.
\end{cor}
\textbf{Proof.} Let $U$ be the subset of $M$ defined by 
$$
U:=\{q\in M\: |\: \widetilde{\Phi}(q)=q\text{ and }T_q\Phi=\text{id}_{T_q\ca}\}.
$$
Remark that if $\widetilde{\Phi}(q)=q$ the equality $T_q\Phi=\text{id}_{T_q\ca}$ is equivalent to the numerical version $T^{\text{num}}_q\Phi=\text{id}_{T^{\text{num}}_q\ca}$. Thus $U$ is a closed subset of $M$.\\

Obviously, the set $U$ is non-empty since $q_0\in U$.  Let $q$ be an element of $U$. By the preceding theorem we have
$$
\displaystyle\Phi\circ\text{exp}\,_q=\text{exp}\,_{\widetilde{\Phi}(q)}\circ T_q\Phi=\text{exp}\,_{q}.
$$
By Theorem \ref{expdiffeo}, the exponential map is a local diffeomorphism and we thus have an open $V$ in $M$ containing $q$ such that $\Phi_{|V}=\text{id}_{\ca_{|V}}$, implying that $U$ is open.
By the connectedness of $M$, it follows that  $U$ equals the whole body $M$ of the supermanifold
$\ca$.\\

Now we know that  $\widetilde{\Phi}(q)=q\text{ and }T_q\Phi=\text{id}_{T_q\ca}$ for all $q\in M$ and consequently we have $\displaystyle\Phi\circ\text{exp}\,_q=\text{exp}\,_{q}$ for all $q\in M$. The morphism $\Phi$ is then
the identity near each $q\in M$. This prove that $\Phi$ is the identity morphism of $\ca$.\hfill$\square$\\

\textbf{Acknowledgments.} The research for this article was partially supported by SFB/TR 12,  ``Symmetries
and Universality in Mesoscopic systems", of the DFG.

\end{document}